\documentclass[12pt]{article}
\usepackage{cite,amsfonts,amsmath,subeqn}
\def\be{\begin{equation}}
\def\eqn#1{\be\label{#1}}

\def\bea{\begin{eqnarray}}
\def\eqnn#1{\bea\label{#1}}
\def\eea{\end{eqnarray}}

\newcommand{\eqna}[1]{\begin{subequations} \label{#1}
\begin{eqnarray}}
\def\eena{\end{eqnarray}
\end{subequations}}
\def\pd{{\partial}} 
\def\hpi{{\hat\pi}} 
\def\({\big(} \def\){\big)}  
\def\y{\eta}\def\vf{\varphi} \def\hv{\hat{\varphi}} 
  
\def\hA{{\hat{A}_{\nu,\r}}} 
\def\hI{{\hat{A}^{-1}_{\nu,\r}}} 

\def\l{\lambda} \def\s{\sigma} \def\r{\rho }\def\t{\tau} 
\def\a{\alpha} \def\b{\beta}
\def\g{\gamma} \def\d{\delta} 
\def\mt{\mapsto}
\def\id{{\bf 1}}
 
\def\spa{~~$\spadesuit$} 
\def\spe{~~\spadesuit} 
\def\PR{{\it Proof: }}

\def\nl{\hfil\break}
\def\nt{\noindent}
\def\lra{\longrightarrow}
\def\Lra{\Longrightarrow}
\def\nn{\nonumber}
\def\cA{{\cal A}}  \def\cC{{\cal C}}
 \def\cE{{\cal E}}

  \def\cU{{\cal U}}

\def\ta{\tilde a}\def\td{\tilde d}
\def\tb{\tilde b}\def\tc{\tilde c}
\def\ha{{\hat a}}\def\hd{{\hat d}}
\def\hb{{\hat b}}\def\hc{{\hat c}}
\def\tA{{\tilde A}} \def\tB{\tilde B}
\def\tC{\tilde C} \def\tD{\tilde D}
\def\hC{{\hat C}}

\def\half{{\textstyle{\frac{1}{2}}}} 
 
\def\ve{\varepsilon}
\def\eps{\epsilon}
\def\om{\omega}

\def\lg{\left\langle}
\def\rg{\right\rangle}

\newcommand{\CC}{\mbox{${\mathbb C}$}}

\newcommand{\ZZ}{\mbox{${\mathbb Z}$}}
\newcommand{\NN}{\mbox{${\mathbb N}$}}

\normalbaselines 
\begin{document}

 \pagestyle{empty}

 \begin{center}
 
 \textsf{\LARGE Duality and Representations for New Exotic Bialgebras}
 
 \vspace{7mm}
 
 {\large D.~Arnaudon$^{a,}$\footnote{Daniel.Arnaudon@lapp.in2p3.fr}, 
 ~A.~Chakrabarti$^{b,}$\footnote{chakra@cpht.polytechnique.fr},\\[2mm] 
 V.K.~Dobrev$^{c,d,}$\footnote{dobrev@inrne.bas.bg,dobrev@ictp.trieste.it} 
 ~and~ S.G.~Mihov$^{c,}$\footnote{smikhov@inrne.bas.bg}}
 
 \vspace{5mm}
 
 \emph{$^a$ Laboratoire d'Annecy-le-Vieux de Physique Th{\'e}orique LAPTH}
 \\
 \emph{CNRS, UMR 5108, associ{\'e}e {\`a} l'Universit{\'e} de Savoie}
 \\
 \emph{LAPTH, BP 110, F-74941 Annecy-le-Vieux Cedex, France}
 \\
 \vspace{3mm}
 \emph{$^b$ Centre de Physique Th{\'e}orique, CNRS UMR 7644}
 \\
 \emph{Ecole Polytechnique, 91128 Palaiseau Cedex, France.}
 \\
 \vspace{3mm}
 \emph{$^c$ Institute of Nuclear Research and Nuclear Energy} 
 \\
 \emph{Bulgarian Academy of Sciences}
 \\
 \emph{72 Tsarigradsko Chaussee, 1784 Sofia, Bulgaria}
 \\
 \vspace{3mm}
 \emph{$^d$ The Abdus Salam International Center for Theoretical Physics}
\\
\emph{Strada Costiera 11, P.O. Box 586}
\\
\emph{34100 Trieste, Italy}
\vspace{3mm}

 \end{center}

 \vspace{.8 cm}
 
\begin{abstract}
We  find the exotic matrix bialgebras which correspond 
to the two non-triangular nonsingular $4\times 4$ $R$-matrices 
in the classification of Hietarinta, namely, ~$R_{S0,3}$~ 
and ~$R_{S1,4}$. We find two new exotic bialgebras ~$S03$~ 
and ~$S14$~ which are not deformations of the 
of the classical algebras of functions on $GL(2)$ or  
$GL(1|1)$. With this we finalize the classification of 
the matrix bialgebras which unital associative algebras 
generated by four elements. We also find the corresponding 
dual bialgebras of these new exotic bialgebras 
and study their representation theory in detail. 
We also discuss in detail a special case of ~$R_{S1,4}$~ 
in which the corresponding algebra turns out to be 
a special case of the two-parameter quantum group deformation 
~$GL_{p,q}(2)$. 
\end{abstract}
 
 \vfill
 
 \rightline{J. Math. Phys. {\bf 43} (2002) 6238-6264;\ {}math.QA/0206053} 
 \rightline{LAPTH-910/02,\ {}RR 036.0502,\ {}INRNE-TH-02-02,\ {}IC/2002/37 
(May 2002)}

 \newpage

\pagestyle{plain}
\setcounter{page}{1}

\section{Introduction}
\label{sect:intro}
\setcounter{equation}{0}

Until now there was no complete list of the matrix bialgebras 
which are unital associative algebras generated by four elements. 
Naturally, since the co-product relations are 
the classical ones we first mention the two related to 
~$GL(2)$, namely, the standard ~$GL_{pq}(2)$\ \cite{DMMZ} 
and nonstandard (Jordanian) $GL_{gh}(2)$ \cite{Ag} two-parameter
deformations. For the supergroup ~$GL(1|1)$~ there are also two: the standard 
$GL_{pq}(1|1)$ \cite{HiRi,DaWa,BuTo} and the hybrid (standard-nonstandard)
$GL_{qh}(1|1)$ \cite{FHR} two-parameter deformations. 
Recently, in \cite{AACDM} it was shown that there are no more 
deformations of $GL(2)$ or $GL(1|1)$.  
In particular, it was shown that these four
deformations match the distinct triangular $4\times 4$
$R$-matrices from the classification of \cite{Hietarinta} 
which are deformations of the trivial $R$-matrix (corresponding
to undeformed $GL(2)$). 

Naturally, there are matrix bialgebras generated by four elements, 
which are not deformations of the classical algebra of functions 
over the group $\ GL(2)$ or the supergroup $GL(1|1)$. 
Those should correspond to $4\times 4$ $R$-matrices 
which are not deformations of the trivial $R$-matrix. 
Studying the classification of \cite{Hietarinta} we noticed altogether 
five nonsingular such $R$-matrices. The triangular ones were introduced in 
\cite{AACDM} and their duals were found and studied in detail in 
\cite{ACDM1}. In the latter paper we called these bialgebras exotic. 

In the present paper we finalize the explicit classification of 
the matrix bialgebras generated by four elements, by studying 
those that correspond to the two  non-triangular nonsingular 
$4\times 4$ $R$-matrices of \cite{Hietarinta}, namely, ~$R_{S0,3}$~ 
and ~$R_{S1,4}$~ which also are not deformations of the trivial $R$-matrix.

The paper is organized as follows. Section 2 just introduces general 
notation. In Section 3 we study the matrix 
bialgebra ~$S03$~ which corresponds to  ~$R_{S0,3}\,$. We find the 
dual bialgebra ~$s03$~ and study the representation theory of  ~$s03$~ 
in detail. In Sections 4 and 5 we study the matrix 
bialgebras ~$S14$~ and ~$S14o$~ which correspond to ~$R_{S1,4}$~ 
for two distinctive regions of the deformation parameter ~$q$~:
~$q^2\neq 1$~ and ~$q^2= 1$, respectively. In both cases we find the 
corresponding dual bialgebras and their representation theory. 
In Section 6 we present our conclusions and outlook. 

\section{Generalities}
\label{sect:gene}
\setcounter{equation}{0}

In this paper we consider  matrix bialgebras 
which are unital associative algebras generated by four elements
~$a,b,c,d$. The co-product and co-unit relations are the classical ones:
\eqna{coal} 
 \delta\;\left( \begin{array}{cc}
 a & b \cr c & d 
 \end{array} \right)
 \ &=&\ \left( \begin{array}{cc}
 a \otimes a + b \otimes c \qquad& a \otimes b + b \otimes d \cr 
 c \otimes a + d \otimes c \qquad& c \otimes b + d \otimes d
 \end{array} \right) \\ 
\ve\ \left( \begin{array}{cc}
a & b \cr c & d 
\end{array} \right)
 \ &=&\ \left(\begin{array}{cc}
1 & 0 \cr 0 & 1
\end{array} \right) 
\eena 
However, the bialgebras under consideration are not Hopf
algebras, except one. This shall be discussed separately in each case. 

It shall be convenient to make the following change of generators: 
\begin{equation}
\label{eq:tatd}
 \tilde{a} \ =\ \half (a+d), \quad
 \tilde{d} \ =\ \half (a-d), \quad
 \tilde{b} \ =\ \half (b+c), \quad
 \tilde{c} \ =\ \half (b-c). 
\end{equation} 
With the new generators we have:
\eqnn{cotilde}
 \delta\;\left( \begin{array}{cc}
 \ta & \tb \cr \tc & \td 
 \end{array} \right)
 \ &=&\ \left( \begin{array}{cc}
 \ta \otimes \ta + \tb \otimes \tb - \tc \otimes \tc + \td \otimes \td 
 \quad& 
 \ta \otimes \tb + \tb \otimes \ta - \tc \otimes \td + \td
 \otimes \tc \cr 
 \ta \otimes \tc + \tc \otimes \ta - \tb \otimes \td + \td
 \otimes \tb 
 \quad& 
 \ta \otimes \td + \td \otimes \ta - \tb \otimes \tc + \tc
 \otimes \tb 
 \end{array} \right) \nn\\
\ve\ \left( \begin{array}{cc}
\ta & \tb \cr \tc & \td 
\end{array} \right)
 \ &=&\ \left(\begin{array}{cc}
1 & 0 \cr 0 & 0
\end{array} \right)
\eea 

\section{Algebra $S03$}
\label{sect:S03}
\setcounter{equation}{0}

\subsection{Bialgebra relations}

In this Section we consider the matrix bialgebra ~$S03$. 
We obtain it by applying the RTT relations of \cite{FRT}: 
\eqn{rtt} R\ T_1\ T_2 \ \ =\ \ T_2\ T_1\ R \ \ , \end{equation} 
where \ $T_1 \ =\ T\, \otimes\, \id_2$\ , \ $T_2 \ =\ \id_2 \,
\otimes\, T$, for the case when $\ R\ =\ R_{S0,3}\ $, where:
\begin{equation}
 \label{eq:S03}
 R_{S0,3}\ \equiv\ 
 \left(
 \begin{array}{cccc}
 1 & 0 & 0 & 1 \cr
 0 & 1 & 1 & 0 \cr
 0 & 1 & -1 & 0 \cr
 -1 & 0 & 0 & 1 \cr 
 \end{array}
 \right)
\end{equation} 
This $R$-matrix is given in \cite{Hietarinta}. 

The relations which follow from (\ref{rtt}) and
(\ref{eq:S03}) are: 
\begin{alignat}{2}
 \label{eq:S03rel}
 & b^2 + c^2 = 0\ , \qquad&&
 a^2 - d^2 = 0\ , \\
 & cd = ba\ , &&
 dc = -ab\ , \nn\\
 & bd = ca\ , &&
 db = -ac\ , \nn\\
 & da = ad\ , &&
 cb = -bc\ . \nn 
\end{alignat}
In terms of the generators ~$\ta,\tb,\tc,\td$~ we have: 
\begin{alignat}{2}
 \label{eq:S03trel}
 & \tb^2 = \tc^2 = 0 \ ,\qquad&&
 \ta\td = \td\ta = 0\ , \\
 & \ta\tb = 0 \ ,&&
 \tb\td = 0\ , \nn\\
 & \td\tc = 0\ , &&
 \tc\ta = 0\ . \nn
\end{alignat}

In view of the above relations we conclude that this bialgebra
has no PBW basis. Indeed, the ordering following from
(\ref{eq:S03trel}) is cyclic: 
\eqn{ord} \ta ~>~ \tc ~>~ \td ~>~ \tb ~>~ \ta 
\end{equation} 
Thus, the basis consists of building blocks like ~$\ta^k\,
\tc\, \td^\ell\, \tb$~ and cyclic. Explicitly the basis can be
described by the following monomials:
\begin{subequations} \label{bas}
\begin{alignat}{2}
&\ta^{k_1}\, \tc\, \td^{\ell_1}\, \tb\, \cdots\, 
 \ta^{k_n}\, \tc\, \td^{\ell_n}\, \tb\, \ta^{k_{n+1}} \ ,&
\qquad & n\,,k_i\,,\ell_i\, \in \, \ZZ_+ \\
&\td^{\ell_1}\, \tb\, \ta^{k_1}\, \tc\, \cdots\, 
 \td^{\ell_n}\, \tb\, \ta^{k_n}\ , && 
n\,, k_i\,,\ell_i\, \in \, \ZZ_+ \\
& \ta^{k_1}\, \tc\, \td^{\ell_1}\, \tb\, \cdots\, 
 \ta^{k_n}\, \tc\, \td^{\ell_n}\ , && 
n\,, k_i\,,\ell_i\, \in \, \ZZ_+ \\
&\td^{\ell_1}\, \tb\, \ta^{k_1}\, \tc\, \cdots\, 
 \td^{\ell_n}\, \tb\, \ta^{k_n}\, \tc\, \td^{\ell_{n+1}} \ , 
&& n\,, k_i\,,\ell_i\, \in \, \ZZ_+ 
\end{alignat}
\end{subequations}

We shall call the elements of the basis 'words'. The one-letter
words are the generators ~$\ta,\tb,\tc,\td$; they are obtained
from ({\ref{bas}}{a}), (\ref{bas}b), (\ref{bas}c), (\ref{bas}d),
resp., for ~$n=0,k_1=1$, ~$n=1, k_1=\ell_1=0$, ~$n=1,
k_1=\ell_1=0$, ~$n=0,\ell_1=1$, resp. The unit element ~$1_\cA$~ is
obtained from (\ref{bas}b) or (\ref{bas}c) for ~$n=0$.

\subsection{Dual algebra}

Two bialgebras \ $\cU , \cA$\ are said to be \ {\it in duality}\ 
\cite{Abe} if there exists a doubly nondegenerate bilinear form 
\eqn{dua} \lg \ , \ \rg \ :\ \cU \times \cA \lra \CC \ , \ \ \ 
\lg \ , \ \rg \ :\ (u,a) \ \mt \ \lg \ u \ , \ a \ \rg
\ , \ u \in\cU \ , \ a \in \cA \end{equation} 
such that, for \ $u,v\in\cU\ , a,b\in\cA$\ :
\eqna{dub} 
&\lg \ u \ , \ ab \ \rg \ =\ \lg \ \d_\cU(u) \ ,
\ a \otimes b \ \rg \ , \ \ \ \lg \ uv \ , \ a \ \rg \ =\ 
\lg \ u \otimes v \ , \ \d_\cA (a) \ \rg \qquad \\ 
&\lg 1_\cU , a \rg \ =\ \ve_\cA (a) \ , \ \ \ \lg u ,
1_\cA \rg \ =\ \ve_\cU (u) \eena \nt 
Two Hopf algebras \ $\cU , \cA$\ are said to be \ {\it in duality}\ 
\cite{Abe} if they are in duality as bialgebras and if 
\be \lg \g_\cU (u) , a \rg \ =\ \lg u , \g_\cA (a)
\rg \end{equation} 

It is enough to define the pairing (\ref{dua}) between the generating
elements of the two algebras. The pairing between any other
elements of \ $\cU, \ \cA$\ follows then from relations (\ref{dub}) and
the standard bilinear form inherited by the tensor product. 

The duality between two bialgebras or Hopf algebras may be used
also to obtain the unknown dual of a known algebra. For that it
is enough to give the pairing between the generating elements of
the unknown algebra with arbitrary elements of the basis of
the known algebra. Using these initial pairings and the duality
properties one may find the unknown algebra. One such possibility
is given in \cite{FRT}. However, their approach is
not universal. In particular, it is not enough for the algebras
considered here, (as will become clear) and will be used only as
consistency check. 

Another approach was initiated by Sudbery \cite{Sudbery}. He
obtained \ $U_q(sl(2)) \otimes U(u(1))$\ as the algebra of
tangent vectors at the identity of \ $GL_q(2)$. The initial pairings
were defined through the tangent vectors at
the identity. However, such calculations become very difficult
for more complicated algebras. Thus, in \cite{Dod} a
generalization was proposed in which the initial pairings are postulated to
be equal to the classical undeformed results. This generalized
method was applied in \cite{Dod} to the standard two-parameter
deformation $GL_{p,q}(2)$, (where also Sudbery's method was
used), then in \cite{DPa} to the multiparameter deformation of
$GL(n)$, in \cite{DPb} to the matrix quantum Lorentz group of
\cite{QL}, in \cite{ADM} to the Jordanian two-parameter
deformation $GL_{g,h}(2)$, in \cite{FHR} to the hybrid two-parameter
deformation of the superalgebra $GL_{q,h}(1|1)$, 
in \cite{DT} to the multiparameter
deformation of the superalgebra $GL(m|n)$, 
in \cite{ACDM1} to the firstly discussed three exotic bialgebras. 
(We note that the dual of $GL_{p,q}(2)$ was obtained also
in \cite{SWZ} by methods of $q$-differential calculus.) 

Let us denote by \ $s03$\ the unknown yet dual algebra of $S03$, 
and by \ $\tA,\tB,\tC,\tD$\ the four generators of \ $s03$. 
Like in \cite{Dod} we define the pairing \ $\lg Z, f\rg$, \
$Z=\tA,\tB,\tC,\tD$, \ $f$\ is from (\ref{bas}), as the
classical tangent vector at the identity: 
\eqn{duda} 
\lg \ Z \ , \ f \ \rg \ \equiv \ 
\ve \left( 
\frac{\partial f}{\partial y} \right) \ , \end{equation} 
where ~$\ (Z,y) \ =\ (\tA,\ta),\ (\tB,\tb),\ (\tC,\tc),\
(\tD,\td)\ $. Explicitly, we get:
\eqna{eq:dualityS03}
 &&\lg \ \tA \ , \ f \ \rg \ =\ 
 \ve\left(\frac{\partial f}{\partial \ta}\right) \ =\ 
 \begin{cases}k \qquad &{\rm for}\ \ f= \ta^k \\ 
 0 \qquad &{\rm otherwise }
 \end{cases} \\ 
 &&\lg \ \tB \ , \ f \ \rg \ =\ 
 \ve\left(\frac{\partial f}{\partial \tb}\right) \ =\ 
 \begin{cases}1 \qquad &{\rm for}\ \ f=\tb \ta^k \\ 
 0 \qquad &{\rm otherwise }
 \end{cases} \\ 
 &&\lg \ \tC \ , \ f \ \rg \ =\ 
 \ve\left(\frac{\partial f}{\partial \tc}\right) \ =\ 
 \begin{cases}1 \qquad &{\rm for}\ \ f= \ta^k \tc \\ 
 0 \qquad &{\rm otherwise }
 \end{cases} \\ 
 &&\lg \ \tD \ , \ f \ \rg \ =\ 
 \ve\left(\frac{\partial f}{\partial \td}\right) \ =\ 
 \begin{cases}
 1 \qquad &{\rm for}\ \ f=\td \\ 
 0 \qquad &{\rm otherwise }
 \end{cases} 
\eena

Using the above we obtain:

\noindent 
{\bf Proposition 1:}\ The generators $\ \tA,\tB,\tC,\tD\ $
introduced above obey the following relations:
\eqna{eq:dualS03}
 && [\tA,Z] = 0\ , \qquad Z=\tB,\tC , \\
& & \tA\tD = \tD\tA = \tD^3 =
 \tB^2 \tD = \tD \tB^2 = \tD , \\
&& [\tB,\tC] = -2\tD , \\
 && \tD\tB = - \tB\tD = \tC\tD^2 =\tD^2 \tC \ , \\
 && \{\tC,\tD\} = 0\ , \\
& & \tB^2 + \tC^2 = 0\ ,\\
 & & \tB^3 = \tB\ , \\
 & & \tC^3 = - \tC\ , \\
& & \tB^2 \tA = \tA\ .
\eena 
\begin{subequations} \label{eq:dualcoS03}
\begin{eqnarray}
 \delta_\cU(\tA) &=& \tA \otimes 1_\cU + 1_\cU\otimes \tA \\
 \delta_\cU(\tB) &=& \tB \otimes 1_\cU + (1_\cU-\tB^2)\otimes \tB \\
 \delta_\cU(\tC) &=& \tC \otimes (1_\cU-\tB^2) + 1_\cU\otimes \tC \\
 \delta_\cU(\tD) &=& \tD \otimes (1_\cU-\tB^2) + (1_\cU-\tB^2)\otimes \tD \\ 
\ve_\cU (Z) &=& 0 \ , \qquad Z \ =\ \tA, \tB, \tC, \tD \ .
\end{eqnarray}
\end{subequations}
\nt $\tA$, ~$\tB^2=-\tC^2$~ and ~$\tD^2$~ are Casimir operators. 
The bialgebra $s03$ is not a Hopf algebra.\nl 
\PR Using the assumed duality the algebraic relations 
(\ref{eq:dualS03}) are shown by calculating their pairings with the 
basis monomials \ $f$\ from (\ref{bas}). In particular, we have 
(giving only the nonzero pairings): 
\eqna{cmmu} 
&&\lg \ \tA \tB \ , \ f \ \rg \ =\ \lg \ \tB \tA \ , \ f \ \rg \ =\ 
k+1\ , \qquad {\rm for}\ \ f=\tb \ta^k\\ 
&&\lg \ \tA \tC \ , \ f \ \rg \ =\ \lg \ \tC \tA \ , \ f \ \rg \ =\ 
k+1\ , \qquad {\rm for}\ \ f=\ta^k \tc\\ 
&&\lg \ \tA \tD \ , \ f \ \rg \ =\ \lg \ \tD \tA \ , \ f \ \rg \ =\ 
\lg \ \tD^3 \ , \ f \ \rg \ = \lg \ \tB^2 \tD \ , \ f \ \rg 
 = \lg \ \tD \tB^2 \ , \ f \ \rg \ =\nn\\ 
&&\ =\ -\lg \ \tB \tC \ , \ f \ \rg \ =\ \lg \ \tC \tB \ , \ f \ \rg \ =\
\lg \ \tD \ , \ f \ \rg \ =\ 1\ , \qquad {\rm for}\ \ f=\td \qquad\\ 
&& \lg \ \tB \tC \ , \ f \ \rg \ =\ \lg \ \tC \tB \ , \ f \
\rg \ =\ 1\ , \qquad {\rm for}\ \ f=\tb\ta^k\tc \\ 
&&\lg \ \tD \tB \ , \ f \ \rg \ =\ -\lg \ \tB \tD \ , \ f \ \rg \ =\ 
\lg \ \tC \tD^2 \ , \ f \ \rg \ =\nn\\ 
&&\ =\ \lg \ \tD^2 \tC \ , \ f \ \rg \ =\ 
1\ , \qquad {\rm for}\ \ f=\tc \\ 
&&\lg \ \tD \tC \ , \ f \ \rg \ =\ -\lg \ \tC \tD \ , \ f \ \rg \ =\ 
 1\ , \qquad {\rm for}\ \ f=\tb \\
&&\lg \ \tB^2 \ , \ f \ \rg \ =\ -\lg \ \tC^2 \ , \ f \ \rg \ =\ 
 1\ , \qquad {\rm for}\ \ f=\ta^k \\
&&\lg \ \tB^3 \ , \ f \ \rg \ =\ \lg \ \tB \ , \ f \ \rg \ =\ 
 1\ , \qquad {\rm for}\ \ f= \tb\ta^k \\
&&\lg \ -\tC^3 \ , \ f \ \rg \ =\ \lg \ \tC \ , \ f \ \rg \ =\ 
 1\ , \qquad {\rm for}\ \ f= \ta^k\tc \\
&&\lg \ \tB^2\tA \ , \ f \ \rg \ =\ \lg \ \tA \ , \ f \ \rg \ =\ 
 k\ , \qquad {\rm for}\ \ f= \ta^k 
\eena 
\nt For the proof of (\ref{cmmu}{h,j}) is used (\ref{cmmu}{g}). The
facts that $\tA$, \ $\tB^2=-\tC^2$\ and \ $\tD^2$\ are Casimir
operators follow easily from (\ref{eq:dualS03}) having in mind also
(\ref{cmmu}{g,j}) and that \ $\langle \tD^2 , \ta \rangle = 1$\ is
the only nonzero pairing of \ $\tD^2$. 
~~For the proof of (\ref{eq:dualcoS03}{a-d}) we use the duality
property (\ref{dub}{a}) 
$$\lg \ Z \ , \ f_1\ f_2 \ \rg \ =\ \lg \ \d_\cU(Z) \ , \ f_1
\otimes f_2 \ \rg $$
for every generator \ $Z$\ of \ $s03$\ and for every \ $f_1, f_2\in
S03$, then we calculate separately the LHS and RHS and compare
the results. We also use that\ 
\eqn{relb}
\tB^2 Z = Z\ , \qquad {\rm for}\ \ Z = \tA,\tB,\tC,\tD
\end{equation} 
- this is contained in (\ref{eq:dualS03}) except for $Z=\tC$ for
which we use that \ $\langle \tB^2\tC , \ta^k\tc \rangle = 1$\ is
the only nonzero pairing of \ $\tB^2\tC$. 
~~Formulae (\ref{eq:dualcoS03}{e}) follow from \
$\ve_\cU (Z) \ =\ \lg \ Z , 1_\cA \ \rg$, cf. (\ref{dub}{b}), and
the defining relations (\ref{eq:dualityS03}). ~~To show that 
$s03$ is not a Hopf algebra we suppose that it is, i.e., there
should exist an antipode $\g$. Then then we use one
of the Hopf algebra axioms:
\eqn{haxi} m \circ ({\rm id} \otimes \g) \circ \d \ =\ i \circ \ve 
\end{equation} 
as maps $s03 \lra s03$, 
where $m$ maps to the usual product in the algebra: $m(Y\otimes Z)
 = YZ\,,\ Y,Z\in s03$ and $i$ is the natural embedding of the
number field $\CC$ into $s03$ : $i(\mu) = \mu 1_\cU\ ,\ \mu\in
\CC$. Applying this to the generator $\tB$ we would get: 
$$ \tB ~+~ \left( 1_\cU -\tB^2 \right) \g (\tB) ~~=~~ 0 $$
which is a contradiction, since\ $1_\cU -\tB^2$\ is zero when
multiplied by anything except by\ $1_\cU$, and in the latter
case the product is not equal to ~$-\tB$. From this we see 
that $\g$ can not be defined on $\tC$ and $\tD$ since their
coproducts involve $\tB$. The antipode may be introduced only 
if we restrict to the subalgebra generated by $\tA$, but the
bialgebra ~$s03$~ as a whole is not a Hopf algebra.\spa

\medskip 

\noindent 
{\bf Corollary 1:}\ The algebra generated by the 
generator $\ \tA\ $ is a sub-bialgebra of ~$s03$. 
The algebra ~$s03'$~ generated by the generators $\ \tB,\tC,\tD\ $ 
is a nine-dimensional sub-bialgebra of ~$s03$~ with PBW basis:
\eqn{pbw:s03} 
1_\cU\,,\ \tB,\ \tC,\ \tD,\ \tB\tC,\ \tB\tD,\ \tD\tC,\ \tB^2,\
\tD^2 \end{equation} 
\PR The statement follows immediately from relations
(\ref{eq:dualS03},\ref{eq:dualcoS03}). We comment only the  
PBW basis of the subalgebra ~$s03'$. Indeed, a priori
it has a PBW basis:
\eqn{pbw:s03a}
\tB^k\,\tD^\ell\,\tC^m \ , \qquad k,\ell\leq 2,\ m\leq 1 
\end{equation} 
the restrictions following from (\ref{eq:dualS03}{b,f,g}).
Furthermore it is easy to see that there are no cubic (and
consequently higher order) elements of the basis. For some of the
cubic elements this is clear from (\ref{eq:dualS03}). For the
rest we have:
\eqna{} 
&&\tB\tD\tC ~=~ -\tD^2\tC^2 ~=~\tD^2\tB^2 ~=~ \tD^2 \\
&&\tB^2\tC ~=~ - \tC^3 ~=~ \tC \\
&&\tB\tD^2 ~=~ - \tC\tD^3 ~=~ \tD\tC \eena 
also using (\ref{eq:dualS03}). Thus, the basis is given by 
(\ref{pbw:s03}) the algebra is indeed nine-dimensional.\spa 

\medskip 

\noindent 
{\bf Remark 1:}\ The algebra ~$s03$~ is not the direct sum of the
two subalgebras described in the preceding Corollary since both
subalgebras have nontrivial action on each other, e.g., 
~$\tB^2\tA ~=~ \tA$, ~$\tA\tD ~=~ \tD$. The algebra ~$s03$~ is a
nine-dimensional associative algebra over the central algebra
generated by $\tA$.\spa 

\subsection{Regular representation} 
\label{subsect:lrr}

We start with the study of the left regular representation (LRR) of the 
subalgebra ~$s03'$. For this we need the left multiplication table:
\begin{center}
 \begin{tabular}{p{1cm}||p{1cm}|p{2.3cm}|p{1.3cm}|p{2.3cm}|p{2.3cm}|p{.7cm}}
   & 1   & $\tB$      & $\tC$   & ${\tB}^2$ &
$\tB\tC$ & $\cdots$ \\[2mm] \hline 
$\Big.\tB$ & $\tB$ & ${\tB}^2$    & $\tB\tC$  & $\tB$   & $\tC$
& $\cdots$ \\[2mm] \hline 
$\Big.\tC$ & $\tC$ & $\tB\tC +2\tD $ & $-{\tB}^2$ & $\tC$   &
$-\tB+2\tD\tC$ & $\cdots$ \\[2mm] \hline 
$\Big.\tD$ & $\tD$ & $-\tB\tD$    & ${\tD}\tC$ & $\tD$   &
$-\tD^2$ & $\cdots$ \\[2mm] \hline 
 \end{tabular}
\end{center}

\vspace{5mm}

\begin{center}
 \begin{tabular}{p{1cm}||p{.7cm}|p{1cm}|p{1.3cm}|p{2.3cm}|p{2.3cm}|}
 & $\cdots$    & $\tD$   & $\tD^2$  & $\tB\tD$ & $\tD\tC$ \\ \hline
 $\Big.\tB$ & $\cdots$ & $\tB\tD$ & $\tD\tC$ & $\tD$   & $\tD^2$ 
 \\[2mm] \hline
 $\Big.\tC$ & $\cdots$ & $-\tD\tC$ & $-\tB\tD$ & $\tD^2$  & $\tD$ 
 \\[2mm] \hline
 $\Big.\tD$ & $\cdots$ & $\tD^2$  & $\tD$   & $-\tD\tC$ & $-\tB\tD$
 \\[2mm] \hline
 \end{tabular}
\end{center}

\vspace{5mm}

The LRR hence contains the subrepresentation generated as a
vector space by ~$\{\tD, \tD^2, \tB\tD, \tD\tC\}$, ~which
decomposes into two two-dimensional irreps 
\begin{align}
 & v^1_0 = \tD+\tD^2\;, && v^1_1 = \tB\tD+\tD\tC\;, \label{eq:irrepS03a}\\
 & \tB \binom{v^1_0}{v^1_1} = \binom{v^1_1}{v^1_0}, &&
 \tC \binom{v^1_0}{v^1_1} = \binom{-v^1_1}{v^1_0}, &&
 \tD \binom{v^1_0}{v^1_1} = \binom{v^1_0}{-v^1_1} \label{eq:irepS03a}\\
{\rm and}&\nn\\ 
 & v^2_0 = \tB\tD-\tD\tC\;, && v^2_1 = \tD - \tD^2 \;, \label{eq:irrepS03b}\\
 & \tB \binom{v^2_0}{v^2_1} = \binom{v^2_1}{v^2_0}, &&
 \tC \binom{v^2_0}{v^2_1} = \binom{-v^2_1}{v^2_0}, &&
 \tD \binom{v^2_0}{v^2_1} = \binom{v^2_0}{-v^2_1} \label{eq:irepS03b}
\end{align}
These two irreps are isomorphic by the map
~$(v^1_0,v^1_1)\to(v^2_0,v^2_1)$. On both of them the Casimirs
~$\tB^2,\tD^2$~ take the value ~1. (Also the Casimir ~$\tA$~ of
$s03$ has the value ~1.) 

The LRR contains also the trivial one-dimensional representation
generated by the vector ~$v ~=~ \tB^2-1_\cU\,$. On this vector
all Casimirs and moreover all generators of $s03$ take the value
~0. 

The quotient of the LRR by the above three sub-modules has the following
multiplication table:
\begin{center}
 \begin{tabular}{p{1cm}||p{1.3cm}|p{1.3cm}|p{1.3cm}|p{1.3cm}|}
    & 1   & $\tB$   & $\tC$   & $\tB\tC$ \\ \hline
 $\Big.\tB$ & $\tB$ & ${\tB}^2 $ & $ \tB\tC $ & $\tC$  \\[2mm] \hline
 $\Big.\tC$ & $\tC$ & $\tB\tC$  & $-{\tB}^2$ & $-\tB$  \\[2mm] \hline
 $\Big.\tD$ & $0$ & $0$ & $0$ & $0$
 \\[2mm] \hline
 \end{tabular}
\end{center}
Thus the quotient decomposes into a direct sum of four one-dimensional
representations, generated as vector spaces by
\begin{equation}
 \label{eq:irrepS03c}
 v_{\eps,\eps'} ~=~ \tB + \eps 1_\cU -i \eps\eps' \tC -i\eps'
\tB\tC \ , \qquad \eps,\eps' = \pm 
\end{equation}
On the latter vectors we have the following action: 
\begin{equation}
 \label{eq:irrepS03d}
 \tB v_{\eps,\eps'} = \eps v_{\eps,\eps'}\ ,  \qquad 
 \tC v_{\eps,\eps'} = i\eps' v_{\eps,\eps'}\ ,  
 \qquad \tD v_{\eps,\eps'} = 0 \;.
\end{equation}
Obviously, on all of them the Casimirs ~$\tB^2,\tD^2$~ take the values
~$1,0$, respectively. However, these four representations are not
isomorphic to each other. 

To summarize, there are seven irreps of ~$s03'$~ 
which are obtained from the LRR:
\begin{itemize}
\item one-dimensional trivial (all generators act by zero)
\item two-dimensional with both Casimirs ~$\tB^2,\tD^2$~
having value 1. 
\item four one-dimensional with Casimir values ~$1,0$~ 
for ~$\tB^2,\tD^2$, respectively. 
\end{itemize}

Turning to the algebra ~$s03$~ we note that it inherits the
representation structure of its subalgebra ~$s03'$. On the
representations (\ref{eq:irrepS03a},\ref{eq:irrepS03b}) the
Casimir ~$\tA$~ has the value 1, while on the trivial irrep
generated by ~$v ~=~ \tB^2-1_\cU$~ the Casimir ~$\tA$~ has the
value 0. However, on the one-dimensional irreps generated by
(\ref{eq:irrepS03c}) the Casimir ~$\tA$~ has no fixed value.
Thus, the list of the irreps of ~$s03$~ arising from the LRR is:
\begin{itemize}
\item one-dimensional trivial 
\item two-dimensional with all Casimirs ~$\tA,\tB^2,\tD^2$~
having value 1. 
\item four one-dimensional with Casimir values ~$\mu,1,0$~ 
for ~$\tA,\tB^2,\tD^2$, respectively, ~$\mu\in \CC$. 
\end{itemize}

 Finally, we note that we could have studied also the right
regular representation of ~$s03$. The list of irreps would be the
same as the one obtained above. 

\subsection{Weight representations} 

Here we consider ~{\it weight representations}. These are
representations which are built from the action of the algebra 
on a ~{\it weight vector}~ with respect to one of the generators.
We start with a weight vector ~$v_0$~ such that:
\eqn{wra} \tD\,v_0 ~=~ \l\,v_0 
\end{equation}
where ~$\l\in \CC$~ is the weight. As we shall see the cases 
$\l\neq 0$ and $\l=0$ are very different.

We start with ~$\l ~\neq ~0$. In that case from from $\ \tD^3
=\tD\ $ follows that ~$\l^2 ~=~ 1$, while from $\ \tB^2\tD =
\tD\ $ follows that ~$\tB^2 v_0 = v_0\,$. Further, from 
(\ref{eq:dualS03}d) follows that ~$\tC\,v_0\ =\ -\l\,\tB\,v_0\,$.  
Thus, acting with the elements
of $s03$ on $v_0$ we obtain a two-dimensional representation, e.g.:
\eqn{wrab} v_0 \ , ~\tB\,v_0 \ , 
\end{equation} 
(and we could have chosen ~$v_0 \ , ~\tC\,v_0$~ as its basis).  
This representation is ~{\it irreducible}. The action is given as
follows:
\begin{center} 
 \begin{tabular}{p{1cm}||p{1.5cm}|p{1.5cm}|}
    & ~$v_0$   & ~$\tB v_0$   \\ \hline
 $\Big.\tB$ & $~\tB v_0$ & ~$v_0 $   \\[2mm] \hline
 $\Big.\tC$ & $-\l \tB v_0$ & ~$\l v_0$  
\\[2mm] \hline 
 $\Big.\tD$ & ~$\l v_0$ & $-\l\tB v_0$ 
 \\[2mm] \hline
 \end{tabular}
\end{center}
Both Casimirs ~$\tB^2,\tD^2$~ take the value ~1. 

Let now ~$\l ~= ~0$. In this case acting with the elements
of $s03$ on $v_0$ we obtain a five-dimensional representation:
\eqn{wrac} v_0 \ , \tB\,v_0 \ ,\tC\,v_0 \ ,\tB\tC\,v_0 \ ,\tB^2\,v_0 \ .
\end{equation} 
This representation is ~{\it reducible}. It has a one-dimensional
subrepresentation spanned by the vector ~$w ~=~ v\, v_0 ~=~ 
(\tB^2-1_\cU)v_0\,$. This is the trivial representation since all
generators act by zero on it. After we factor out this 
representation the factor-representation splits into four 
one-dimensional representations spanned by the following vectors
~$w_{\eps,\eps'} ~=~ v_{\eps,\eps'}\,v_0\,$, where
~$v_{\eps,\eps'}$~ are from (\ref{eq:irrepS03c}) and the action
of the generators is as given in (\ref{eq:irrepS03d}). Thus,
these irreps are as those obtained from the LRR. 

To summarize, there are six irreps of ~$s03'$~ 
which are obtained as weight irreps of the generator ~$\tD$~:
\begin{itemize}
\item one-dimensional trivial 
\item one two-dimensional with both Casimirs ~$\tB^2,\tD^2$~
having value 1. 
\item four one-dimensional with Casimir values ~$1,0$~ 
for ~$\tB^2,\tD^2$, respectively. 
\end{itemize}

Turning to the algebra ~$s03$~ we note that it inherits the
representation structure of its subalgebra ~$s03'$, however, the
value of the Casimir ~$\tA$~ is not fixed except on the trivial
irrep. Thus, the list of the irreps of ~$s03$~ 
which are obtained as weight irreps of the generator ~$\tD$~ is:
\begin{itemize}
\item one-dimensional trivial 
\item one two-dimensional with Casimir values ~$\mu,1,1$~ 
for ~$\tA, \tB^2,\tD^2$, respectively, $\mu\in \CC$. 
\item four one-dimensional with Casimir values ~$\mu,1,0$~ 
for ~$\tA, \tB^2,\tD^2$, respectively, $\mu\in \CC$. 
\end{itemize}

Finally, we note that it is not possible to construct weight
representations w.r.t. generator ~$\tB$ (or $\tC$). 

\subsection{Representations of s03 on S03}
\label{ind:s03}

Here we shall study the representations of ~$s03$~ obtained by
the use of its right regular action (RRA) on the dual bialgebra ~$S03$.
 The RRA is defined as follows:
\eqna{rraa} 
&& \pi_R(Z) f ~\equiv~ f_{(1)} \lg Z,f_{(2)}\rg \ , \qquad Z\in
s03\,, \ Z\neq 1_\cU \,, \quad f\in S03 \ , \\
&&\pi_R(1_\cU) f ~\equiv~ f \ , \qquad f\in S03 \ ,
\eena 
where we use Sweedler's notation for the co-product: ~$\d(f) ~=~
f_{(1)} \otimes f_{(2)}\,$. (Note that we can not use the left
regular action since that would be given by the formula:
~$\pi_L(Z) f = \lg \g_\cU (Z),f_{(1)}\rg\,f_{(2)}$~ and 
we do not have an antipode.) More explicitly, for the generators
of ~$s03$~ we have: 
\eqna{rrag}
 &&
 \pi_R(\tA) 
 \begin{pmatrix} 
 \ta & \tb \cr \tc & \td 
 \end{pmatrix}
 = 
 \begin{pmatrix} 
 \ta & \tb \cr \tc & \td 
 \end{pmatrix}
 \\
&& \pi_R(\tB) 
 \begin{pmatrix} 
 \ta & \tb \cr \tc & \td 
 \end{pmatrix}
 = 
 \begin{pmatrix} 
 \tb & \ta \cr \td & \tc 
 \end{pmatrix}
 \\ &&
 \pi_R(\tC) 
 \begin{pmatrix} 
 \ta & \tb \cr \tc & \td 
 \end{pmatrix}
 = 
 \begin{pmatrix} 
 -\tc & \td \cr \ta & -\tb 
 \end{pmatrix}
 \\ &&
 \pi_R(\tD) 
 \begin{pmatrix} 
 \ta & \tb \cr \tc & \td 
 \end{pmatrix}
 = 
 \begin{pmatrix} 
 \td & -\tc \cr -\tb & \ta 
 \end{pmatrix} \\
 && \pi_R(Z)\,1_\cA ~=~ 1_\cA\, \lg Z,1_\cA\rg ~=~ 1_\cA\,
\ve_\cU (Z) ~=~ 0\ , \qquad Z=\tA,\tB,\tC,\tD \qquad 
\eena 

For the action on the elements (words) of ~$S03$ we use a
Corollary of (\ref{rraa}):
\eqn{rrap} 
 \pi_R(Z) fg ~=~ \pi_R(\d_\cU(Z)) (f\otimes g) 
\end{equation} 
where ~$f,g$~ are arbitrary words from (\ref{bas}). 
Further we shall need the notion of the 'length' $\ell(f)$ of the
word $f$. It is defined naturally as the number of the letters of
$f$; in addition we set ~$\ell (1_\cA) ~=~ 0$. 
Now we obtain from (\ref{rrap}):
\eqna{rrapp}
&& \pi_R(\tA)\, f ~=~ \ell (f)\, f \\
&& \pi_R(\tB)\, f\cdot g ~=~ (\pi_R(\tB) f)\cdot g \\
&& \pi_R(\tC)\, f\cdot g ~=~ f\cdot (\pi_R(\tC) g) \\
&& \pi_R(\tD)\, f ~=~ 0 \ , \qquad {\rm if}\ \ell (f) > 1 
\eena 
{}From (\ref{rrapp}{b,c}) it is obvious that the only nonzero action
of ~$\tB,\tC$~ actually is: 
\eqna{rraw}
 && \pi_R(\tB) \begin{pmatrix} 
 \ta & \tb \cr \tc & \td 
 \end{pmatrix} \cdot f 
 = 
 \begin{pmatrix} 
 \tb & \ta \cr \td & \tc 
 \end{pmatrix}
 \cdot f
 \\
 && \pi_R(\tC)\, f\cdot
 \begin{pmatrix} 
 \ta & \tb \cr \tc & \td 
 \end{pmatrix}
 = 
f\cdot
 \begin{pmatrix} 
 -\tc & \td \cr \ta & -\tb 
 \end{pmatrix} 
\eena
 
{}From (\ref{rrapp}{a}) it is obvious that we can classify the
irreps by the value ~$\mu_A$~ of the Casimir ~$\tA$~ which runs over the
nonnegative integers. For fixed ~$\mu_A$~ the basis of the
corresponding representations is spanned by the words ~$f$~ such
that ~$\ell (f) ~=~ \mu_A\,$. Thus, we have:

\begin{itemize}
\item $\mu_A ~=~ 0$ \\
This is the one-dimensional trivial representation spanned by the
unit element ~$1_\cA$~ on which all generators of ~$s03$~ have
zero action. 
\item $\mu_A ~=~ 1$ \\ 
This representation is four-dimensional spanned by the four
generators ~$\ta,\tb,\tc,\td$~ of ~$S03$. It is reducible and
decomposes in two two-dimensional irreps with basis vectors: 
\eqna{eq:irrepS03aa}
& v^1_0 ~=~ \ta + \td ~=~ a\ , \quad v^1_1 ~=~ \tb + \tc ~=~ b\ , \\
{\rm and}&\nn\\ 
& v^2_0~=~ \tb - \tc ~=~ c \ , \quad 
v^2_1 ~=~ \ta - \td ~=~ d\ . 
\eena 
The RRA of ~$\tB,\tC,\tD$~ on these vectors is as 
(\ref{eq:irepS03a},\ref{eq:irepS03b})~:  
\eqn{eq:irepS03} 
 \pi_R(\tB) \binom{v^k_0}{v^k_1} = \binom{v^k_1}{v^k_0}, \qquad
\pi_R(\tC) \binom{v^k_0}{v^k_1} = \binom{-v^k_1}{v^k_0}, \qquad
\pi_R(\tD) \binom{v^k_0}{v^k_1} = \binom{v^k_0}{-v^k_1} 
\end{equation}
These two irreps are isomorphic by the map
~$(v^1_0,v^1_1)\to(v^2_0,v^2_1)$. On both of them the Casimirs
~$\tB^2,\tD^2$~ take the value ~1. 

\item $\mu_A ~=~ 2$ \\ 
This representation is eight-dimensional spanned by ~$\ta^2,\ta\tc,
\tb\ta,\tb\tc,\tc\tb,\tc\td,\td^2,\td\tb$. It is reducible and
decomposes in eight one-dimensional irreps with basis vectors: 
\eqna{eq:onedima}
& v^1_{\eps,\eps'} ~=~ (\ta +\eps\tb)(\ta +i\eps'\tc) \\ 
& v^2_{\eps,\eps'} ~=~ (\td +\eps\tc)(\td +i\eps'\tb) \\ 
&\eps,\eps' ~=~ \pm \nn 
\eena 
The RRA of ~$\tB,\tC,\tD$~ on these vectors is as 
 (\ref{eq:irrepS03d})~: 
\eqn{onedim} 
 \pi_R(\tB) v^k_{\eps,\eps'} = \eps v^k_{\eps,\eps'}\ ,  \qquad 
 \pi_R(\tC) v^k_{\eps,\eps'} = i\eps' v^k_{\eps,\eps'}\ ,  
 \qquad \pi_R(\tD) v^k_{\eps,\eps'} = 0 \;.
\end{equation}
The irrep with vector ~$v^1_{\eps,\eps'}$~ is isomorphic to 
the irrep with vector ~$v^2_{\eps,\eps'}\,$. Thus, there are only
four distinct irreps parametrized by ~$\eps,\eps'$. On all of them
the Casimirs ~$\tB^2,\tD^2$~ take the value ~1,0, respectively. 

\item $\mu_A ~=~ N > 2$ \\ 
These representations are reducible and
decompose in one-dimensional irreps with basis vectors: 
\eqna{eq:ondima}
& v^1_{\eps,\eps'} ~=~ (\ta +\eps\tb)\cdot f_1 \cdot(\ta +i\eps'\tc) \\ 
 \label{eq:ondimb}
& v^2_{\eps,\eps'} ~=~ (\td +\eps\tc)\cdot f_2 \cdot(\td +i\eps'\tb) \\ 
 \label{eq:ondimc}
& v^3_{\eps,\eps'} ~=~ (\ta +\eps\tb)\cdot f_3 \cdot(\td +i\eps'\tb) \\ 
 \label{eq:ondimd}
& v^4_{\eps,\eps'} ~=~ (\td +\eps\tc)\cdot f_4 \cdot(\ta +i\eps'\tc) \\ 
&\eps,\eps' ~=~ \pm \ , \qquad \ell (f_k) = N-2 \nn 
\eena 
The RRA of ~$\tB,\tC,\tD$~ on these vectors is as exactly as
(\ref{eq:irrepS03d}). The irrep with vector ~$v^k_{\eps,\eps'}$~ is
isomorphic to the irrep with vector ~$v^n_{\eps,\eps'}\,$. Thus,
there are only four distinct irreps as in the case above. 
On all of them the Casimirs ~$\tB^2,\tD^2$~ take the value
~1,0, respectively. 
\end{itemize}

To summarize the list of irreps of ~$s03'$~ is the same as given
in subsection \ref{subsect:lrr}. The list of irreps of ~$s03$~ here is
smaller since the Casimir ~$\tA$~ can take only nonnegative integer values.
Thus, the list of the irreps of ~$s03$~ using the dual bialgebra
~$S03$~ as carrier space is:
\begin{itemize}
\item one-dimensional trivial
\item two-dimensional with all Casimirs ~$\tA,\tB^2,\tD^2$~
having value 1. 
\item four one-dimensional with Casimir values ~$\mu,1,0$~ 
for ~$\tA,\tB^2,\tD^2$, respectively, $\mu\in\NN+1$. 
\end{itemize}
The difference in the two lists is natural since here 
more structure (the co-product) is involved. Speaking more loosely
the irreps here may be looked upon as 'integrals' of the irreps
obtained in subsection \ref{subsect:lrr}.

\section{Algebra $S14$}
\label{sect:S14}
\setcounter{equation}{0}

\subsection{Bialgebra relations}

In this Section we consider the matrix bialgebra ~$S14$. 
We obtain it by applying the RTT relations (\ref{rtt})
for the case $\ R\ =\ R_{S1,4}\ $, when ~$q^2\neq 1$~ where:
\begin{equation}
\label{eq:S14}
 R_{S1,4}\ \equiv\ \left(
 \begin{array}{cccc}
 0 & 0 & 0 & q \cr
 0 & 0 & 1 & 0 \cr
 0 & 1 & 0 & 0 \cr
 q & 0 & 0 & 0 \cr 
 \end{array}
 \right)
\end{equation} 
This $R$-matrix is given in \cite{Hietarinta}. 

The relations which follow from (\ref{rtt}) and
(\ref{eq:S14}) when $q^2 \ne 1$ are: 
\begin{align}
 \label{eq:S14rel}
 & b^2 - c^2 = 0\ , &&
 a^2 - d^2 = 0 \\
 & ab = ba = 0\ , &&
 ac = ca = 0 \nn\\
 & bd = db = 0\ , &&
 cd = dc =0 \nn 
\end{align}
In terms of the generators $\ta,\tb,\tc,\td$ 
\begin{align}
 \label{eq:S14trel}
 & \tb\tc + \tc\tb = 0 &&
 \ta\td + \td\ta = 0 \\
 & \ta\tb = \tb\ta = 0 &&
 \ta\tc = \tc\ta = 0 \nn\\
 & \tb\td = \td\tb = 0 &&
 \tc\td = \td\tc = 0 \nn 
\end{align}
{}From the above relations it is clear that the PBW basis of
$S14$ is:
\eqn{bas14} \ta^k\td^\ell \ , \quad \tb^k\tc^\ell 
\end{equation}

\subsection{Dual algebra}

Let us denote by \ $s14$\ the unknown yet dual algebra of $S14$, 
and by ~$\tA,\tB,\tC,\tD$~ the four generators of \ $s14$. 
We define the pairing \ $\lg Z, f\rg$, \
$Z=\tA,\tB,\tC,\tD$, \ $f$\ is from (\ref{bas14}), as
(\ref{duda}). Explicitly, we obtain:
\eqna{eq:dualityS14}
 &&\lg \ \tA \ , \ f \ \rg \ =\ 
 \ve\left(\frac{\partial f}{\partial \ta}\right) \ =\ 
 \begin{cases}k\delta_{\ell 0} \qquad &f= \ta^k \td^\ell\\ 
 0 \qquad &f=\tb^k \tc^\ell 
 \end{cases} \\ 
 &&\lg \ \tB \ , \ f \ \rg \ =\ 
 \ve\left(\frac{\partial f}{\partial \tb}\right) \ =\ 
 \begin{cases}
0 &f= \ta^k \td^\ell\\ 
\delta_{k1}\delta_{\ell 0} \qquad &f=\tb^k \tc^\ell 
 \end{cases} \\ 
 &&\lg \ \tC \ , \ f \ \rg \ =\ 
 \ve\left(\frac{\partial f}{\partial \tc}\right) \ =\ 
 \begin{cases}
 0&f= \ta^k \td^\ell\\ 
 \delta_{k0}\delta_{\ell 1} \qquad &f=\tb^k \tc^\ell \\ 
 \end{cases} \\ 
 &&\lg \ \tD \ , \ f \ \rg \ =\ 
 \ve\left(\frac{\partial f}{\partial \td}\right) \ =\ 
 \begin{cases}
 \delta_{\ell 1} \qquad &f=\ta^k \td^\ell \\ 
 0 \qquad &f=\tb^k \tc^\ell 
 \end{cases} 
\eena 

We shall need (as in \cite{ACDM1}) the auxilliary operator ~$E$~
such that
\eqn{axi} 
\left\langle E, f \right\rangle \ =\ 
\begin{cases}1 &{\rm for}\ \ f=1_\cA\\ 
0 &{\rm otherwise} 
\end{cases} 
\end{equation}

Using the above we obtain:

\medskip

\noindent 
{\bf Proposition 2:}\ The generators $\ \tA,\tB,\tC,\tD\ $
introduced above obey the following relations:
\eqna{eq:dualS14}
&& \tC = \tD\tB = -\tB\tD \\ 
&& [\tA,\tD] = 0 \\
&& \tA\tB = \tB\tA = \tD^2 \tB = \tB^3 = \tB \\
&& E Z \ =\ Z E \ =\ 0 \ , \quad Z=\tA,\tB,\tD \ . 
\eena 
\eqna{eq:dualco14}
 \delta_\cU(\tA) &=& \tA \otimes 1_\cU + 1_\cU\otimes \tA \\
 \delta_\cU(\tB) &=& \tB \otimes E + E\otimes \tB \\
 \delta_\cU(\tD) &=& \tD \otimes K + 1_\cU\otimes \tD \ , \qquad 
 K \equiv (-1)^{\tA} \\ 
 \delta(E) &=& E\otimes E \eena
 \eqna{cougq} 
\ve_\cU (Z) &=& 0 \ , \qquad Z \ =\ \tA, \tB, \tD \\ 
\ve_\cU (E) &=& 1
\eena 
$\tA$, $\tB^2$ and $\tD^2$ are Casimir operators. The bialgebra
$s14$ is not a Hopf algebra.\nl 
\PR Using the assumed duality the algebraic relations 
(\ref{eq:dualS14}) are shown by calculating their pairings with the 
basis monomials \ $f$\ from (\ref{bas14}). In particular, we have 
(giving only the nonzero pairings): 
\eqna{cm14} 
&&\lg \ \tD \tB \ , \ f \ \rg \ =\ -\lg \ \tB \tD \ , \ f \ \rg \ =\ 
\lg \ \tC \ , \ f \ \rg \ =\ 1\ , \qquad {\rm for}\ \ f=\tc\quad \\ 
&&\lg \ \tA \tD \ , \ f \ \rg \ =\ \lg \ \tD \tA \ , \ f \ \rg \
=\ k+1\ , \qquad {\rm for}\ \ f=\ta^k \td \\ 
&&\lg \ \tA \tB \ , \ f \ \rg \ =\ \lg \ \tB \tA \ , \ f \ \rg \ =\ 
\lg \ \tD^2\tB \ , \ f \ \rg \ =\nn\\ &&=\ \lg \ \tB^3 \ , \ f \ \rg \ =\ 
\lg \ \tB \ , \ f \ \rg \ =\ 1\ , \qquad {\rm for}\ \ f=\tb 
\eena
The facts that $\tA$, \ $\tB^2$\ and \ $\tD^2$\ are Casimir
operators follow easily from (\ref{eq:dualS14}).
The proof of (\ref{eq:dualco14}) is done as the proof of 
(\ref{eq:dualcoS03}). We also use that 
\eqn{reln}
\lg \ \tA^n \ , \ \ta^k \td^\ell\ \rg \ =\ k^n \d_{\ell 0} 
\end{equation}
and hence:
\eqn{rel1}
\lg \ K \ , \ \ta^k \td^\ell\ \rg \ =\ (-1)^k \d_{\ell 0} 
\end{equation}
There is no antipode for the bialgebra $s14$. Indeed, suppose
that there was such. Then by applying the Hopf algebra axiom
(\ref{haxi}) to the operator $E$ we would get:
$$ E\, \g(E)\ =\ 1_\cU $$ 
which would lead to contradiction after multiplication from the
left with $\ Z \ =\ \tA, \tB, \tD\ $ (we would get $0=Z$). From
this follows also that the generator ~$\tB$~ does not have an
antipode, since from (\ref{haxi}) to the $\tB$ we would get:
$$ \tB\, \g(E)\ + E\, \g(\tB) =0 $$ 
Thus, the bialgebra ~$s14$~ is not a Hopf algebra.\spa

\medskip 

\noindent 
{\bf Corollary 2:}\ The algebra generated by the 
generator $\ \tA\ $ is a sub-bialgebra of ~$s14$. 
The algebra ~$s14'$~ generated by $\ \tB,\tD\ $ 
is a subalgebra of ~$s14$, but is not a sub-bialgebra 
(cf. (\ref{eq:dualco14}{b,c}). It has the following PBW basis:
\eqn{pbw:s14} 
\tB,\ \tB^2,\ \tD\tB,\ \tD\tB^2,\ \tD^\ell,\ \ell = 0,1,2,...
 \end{equation} 
where we use the convention ~$\tD^0 = 1_\cU\,$.\spa 

\subsection{Regular representation} 

We start with the study of the right regular representation of
the subalgebra ~$s14'$. For this we use the right multiplication table: 
\begin{center}
\begin{tabular}{p{1cm}||p{1cm}|p{2.3cm}|p{1.3cm}|p{2.3cm}|p{2.3cm}|p{.7cm}}
   & $\tB$   & $\tB^2$  & $\tD\tB$  & $\tD{\tB}^2$ &
$\tD^{2k}$  & $\tD^{2k+1}$  \\[2mm] \hline 
$\Big.\tB$ & $\tB^2$  & $\tB$   & $\tD\tB^2$ & $\tD\tB$   & 
$\tB$    & $\tD\tB$   \\[2mm] \hline 
$\Big.\tD$ & $-\tD\tB$ & $\tD\tB^2$ & $-\tB$   & $\tB^2$ &
$\tD^{2k+1}$ & $\tD^{2k+2}$  \\[2mm] \hline 
 \end{tabular}
\end{center}

{}From the above table follows that there is a four-dimensional
subspace spanned by ~$\tB, \tB^2, \tD\tB, \tD\tB^2$. It is
reducible and decomposes into four one-dimensional
representations spanned by:
\eqn{oned14} 
v_{\eps,\eps'} ~=~ \tB + \eps \tB^2 - \eps' \tD\tB + \eps\eps'
\tD\tB^2 
\end{equation} 
The action of ~$\tB,\tD$~ on these vectors is:
\eqn{act14} \tB v_{\eps,\eps'} ~=~ \eps v_{\eps,\eps'} \ , \qquad
\tD v_{\eps,\eps'} ~=~ \eps' v_{\eps,\eps'} 
\end{equation} 
The value of the Casimirs ~$\tB^2,\tD^2$~ on these vectors is ~1. 

The quotient of the RRR by the above submodules has the following
multiplication table:
\begin{center}
 \begin{tabular}{p{1cm}||p{1.3cm}|p{1.3cm}|p{1.3cm}|p{1.3cm}|}
    & $\tD^{2k}$ & $\tD^{2k+1}$  \\[2mm] \hline 
 $\Big.\tB$ & 0     & 0   \\[2mm] \hline
 $\Big.\tD$ & $\tD^{2k+1}$ & $\tD^{2k+2}$  \\[2mm] \hline 
 \end{tabular}
\end{center}
This representation is reducible. It contains an infinite set of
nested submodules ~$V^n\supset V^{n+1}$, $n=0,1,...$, where
~$V^n$~ is spanned by ~$\tD^{n+\ell}$, $\ell=0,1,...$.
Correspondingly there is an infinite set of one-dimensional
irreducible factor-modules ~$F^n\equiv V^n/V^{n+1}$, (generated
by ~$\tD^n$) which are
all isomorphic to the trivial representation since the generators
~$\tB,\tD$~ act as zero on them. Thus there are five irreps
arising from the RRR of ~$s14'$~:
\begin{itemize}
\item one-dimensional trivial 
\item four one-dimensional with both Casimirs ~$\tB^2,\tD^2$~
having value 1. 
\end{itemize}

Turning to the algebra ~$s14$~ we note that it inherits the
representation structure of its subalgebra ~$s14'$. On the 
representations (\ref{oned14}) the
Casimir ~$\tA$~ has the value 1. However, on the one-dimensional
irreps ~$F^n$~ the Casimir ~$\tA$~ has no fixed value.
Thus, the list of the irreps arising from the RRR of ~$s14$~ is: 
\begin{itemize}
\item one-dimensional with Casimir values ~$\mu,0,0$~ 
for ~$\tA,\tB^2,\tD^2$, respectively, ~$\mu\in \CC$. 
\item four one-dimensional with all Casimirs ~$\tA,\tB^2,\tD^2$~
having value 1. 
\end{itemize}

\subsection{Weight representations} 

Here we study weight representations, first w.r.t. $\tD$, as in (\ref{wra}).
The resulting representation of ~$s14'$~ is three-dimensional: 
\eqn{wr14} v_0 \ , ~\tB\,v_0 \ , ~\tB^2\,v_0 \ .
\end{equation} 
It is reducible and contains one one-dimensional and one
two-dimensional irrep:
\begin{itemize}
\item one-dimensional
\begin{align}
 \label{eq:irrepS14a}
 & w_0 ~=~ (\tB^2-1_\cU)v_0\;, && \\
&\tB w_0 =0 \,, && \tD {w_0} = \lambda w_0\,, 
\end{align}
$\lambda\in \CC$. 
\item two-dimensional
\begin{align}
 \label{eq:irrepS14b}
 & \{v_0,v_1=\tB\,v_0\} &&\\[2mm]
 & \tB \binom{v_0}{v_1} = \binom{v_1}{v_0}, &&
 \tD \binom{v_0}{v_1} = \lambda \binom{v_0}{-v_1}
\end{align}
with $\lambda=\pm 1$.
\end{itemize}

Turning to the algebra ~$s14$~ we note that it inherits the
representation structure of its subalgebra ~$s14'$. On the
one-dimensional irrep (\ref{eq:irrepS14a}) the Casimir ~$\tA$~
has no fixed value since $\tB$ is trivial, and $[\tA,\tD]=0$. 
On the two-dimensional irrep (\ref{eq:irrepS14b}) 
the Casimir ~$\tA$~ has the value 1 since ~$\tA\tB = \tB$. 

Thus, there are the following irreps of ~$s14$~ 
which are obtained as weight irreps of the generator ~$\tD$~:
\begin{itemize}
\item one-dimensional with Casimir values ~$\mu,0,\l^2$~ 
for ~$\tA, \tB^2,\tD^2$, respectively, $\mu,\l\in \CC$. 
\item two two-dimensional with all Casimirs 
~$\tA, \tB^2,\tD^2$~ having the value 1. 
\end{itemize}

Next we consider weight representations w.r.t. $\tB$~:
\eqn{wr14b} \tB\,v_0 ~=~ \nu\,v_0 \ ,
\end{equation}
with $\nu\in \CC$. From $\tB^3 = \tB$ follows that ~$\nu=0,\pm 1$. 
Acting with the generators we obtain the
following representation vectors: ~$v_\ell ~=~ \tD^\ell\,v_0\,$. 
We have that ~$\tD v_\ell ~=~ v_{\ell+1}$. 

Further we consider first the case ~$\nu^2=1$. Then we apply the
relation ~$\tD^2 \tB = \tB$ to $v_\ell$ and we get: 
$$\tD^2 \tB v_\ell ~=~ (-1)^\ell \nu v_{\ell+2} ~=~ 
\tB v_\ell ~=~ (-1)^\ell \nu v_{\ell} $$ 
from which follows that we have to identify $v_{\ell+2}$ with
$v_\ell\,$. Thus the representation is given as follows:
\begin{align}
 \label{eq:irrepS14bb}
 & \{v_0,v_1=\tD\,v_0\} &&\\[2mm]
 & \tB \binom{v_0}{v_1} = \nu \binom{v_0}{-v_1}, &&
 \tD \binom{v_0}{v_1} = \binom{v_1}{v_0}
\end{align}
On this irrep both Casimirs ~$\tB^2,\tD^2$~ have value ~1,
($\nu^2=1$). 

Further we consider the case ~$\nu=0$. This representation is
reducible. It contains an infinite set of 
nested submodules ~$V^n\supset V^{n+1}$, $n=0,1,...$, where
~$V^n$~ is spanned by ~$\tD^{n+\ell} v_0$, $\ell=0,1,...$.
Correspondingly there is an infinite set of one-dimensional
irreducible factor-modules ~$F^n\equiv V^n/V^{n+1}$, (generated
by ~$\tD^n v_0$) which are
all isomorphic to the trivial representation since the generators
~$\tB,\tD$~ act as zero on them.

Turning to the algebra ~$s14$~ we note that it inherits the
representation structure of its subalgebra ~$s14'$, with the
value of the Casimir ~$\tA$~ being not fixed if $\tB$ acts
trivially, and being 1, if $\tB$ acts non trivially.

Thus, there are the following irreps of ~$s14$~ 
which are obtained as weight irreps of the generator ~$\tB$~:
\begin{itemize}
\item one-dimensional with Casimir values ~$\mu,0,0$~ 
for ~$\tA, \tB^2,\tD^2$, respectively, $\mu\in \CC$. 
\item two two-dimensional with all Casimirs 
~$\tA, \tB^2,\tD^2$~ having the value 1.
\end{itemize}

\subsection{Representations of s14 on S14}
\label{ind:s14}

Here we shall study the representations of ~$s14$~ obtained by
the use of its right regular action (RRA) on the dual bialgebra ~$S14$.
 The RRA is defined as in (\ref{rraa}). For the generators
of ~$s14$~ we have:  
\eqna{rra14}
 &&
 \pi_R(\tA) 
 \begin{pmatrix} 
 \ta & \tb \cr \tc & \td 
 \end{pmatrix}
 = 
 \begin{pmatrix} 
 \ta & \tb \cr \tc & \td 
 \end{pmatrix}
 \\
&& \pi_R(\tB) 
 \begin{pmatrix} 
 \ta & \tb \cr \tc & \td 
 \end{pmatrix}
 = 
 \begin{pmatrix} 
 \tb & \ta \cr \td & \tc 
 \end{pmatrix} \\ 
 && \pi_R(\tD) 
 \begin{pmatrix} 
 \ta & \tb \cr \tc & \td 
 \end{pmatrix}
 = 
 \begin{pmatrix} 
 \td & -\tc \cr -\tb & \ta 
 \end{pmatrix} \\
 && \pi_R(E) 
 \begin{pmatrix} 
 \ta & \tb \cr \tc & \td 
 \end{pmatrix}
 = 
 \begin{pmatrix} 
 0 & 0 \cr 0 & 0 
 \end{pmatrix} \\
 && \pi_R(Z)\,1_\cA ~=~ 1_\cA\, \lg Z,1_\cA\rg ~=~ 1_\cA\,
\ve_\cU (Z) ~=~ 
\begin{cases} 
0\ , \quad & Z=\tA,\tB,\tD \\ 
1\ , \quad & Z=E
\end{cases} 
\eena 
For the action on the basis of ~$S14$ we use formula
(\ref{rrap}). We obtain: 
\eqna{act14r}
 && \pi_R(A)\, \ta^n \td^k = (n+k) \ta^{n} \td^{k}\ ,
 \qquad \pi_R(A)\, \tb^n \tc^k = (n+k) \tb^{n} \tc^{k}
 \\ 
 && \pi_R(B)\, \ta^n \td^k = \delta_{k0}\delta_{n1} \tb
 + \delta_{n0}\delta_{k1} \tc \ ,
 \qquad \pi_R(B)\, \tb^n \tc^k = \delta_{k0}\delta_{n1} \ta
 + \delta_{n0}\delta_{k1} \td
\qquad \\ 
 && \pi_R(D)\, \ta^k \td^\ell = (-1)^{\ell+1}\ell\ta^{k+1}\td^{\ell-1} + 
(-1)^{\ell}k\ta^{k-1}\td^{\ell+1}\\ 
 && \pi_R(D)\, \tb^k \tc^\ell = (-1)^{\ell}\ell\tb^{k+1}\tc^{\ell-1} + 
(-1)^{\ell+1}k\tb^{k-1}\tc^{\ell+1} 
\eena
We see that similarly to subsection \ref{ind:s03} the Casimir
$\tA$ acts as the length of the elements of $S14$, i.e.,
(\ref{rrag}) holds. Thus, also here we classify the 
irreps by the value ~$\mu_A$~ of the Casimir ~$\tA$~ which runs over the
nonnegative integers. For fixed ~$\mu_A$~ the basis of the
corresponding representations is spanned by the elements ~$f$~ such
that ~$\ell (f) ~=~ \mu_A\,$. The dimension of each such
representation is:
\eqn{dim14} {\rm dim}\, (\mu_A) ~=~ 
\begin{cases} 
2 (\mu_A +1) \qquad &{\rm for}\ \ \mu_A\geq 1 \\ 
1 \qquad &{\rm for}\ \ \mu_A =0 
 \end{cases}
\end{equation}

The classification goes as follows: 
\begin{itemize}
\item $\mu_A ~=~ 0$ \\ 
This is the one-dimensional trivial representation spanned by
~$1_\cA\,$. 
\item $\mu_A ~=~ 1$ \\ 
This representation is four-dimensional spanned by the four 
generators ~$\ta,\tb,\tc,\td$~ of ~$S14$. It decomposes in two
two-dimensional isomorphic to each other irreps with basis
vectors as in (\ref{eq:irrepS03aa}) - this is due to the fact
that the action (\ref{rra14}{b,c}) is the same as the action
(\ref{rrag}). The value of the Casimirs ~$\tB^2,\tD^2$~ is 1. 
\item Each representation for fixed ~$\mu_A\geq 2$~ is reducible and
decomposes in two isomorphic representations: one built on the basis
~$\ta^k\td^\ell$, and the other built on the basis
~$\tb^k\tc^\ell$, each of dimension ~$\mu_A +1$. Thus, for
~$\mu_A\geq 2$~ we shall consider only the
representations built on the basis ~$\ta^k\td^\ell$. 
These representations are also reducible and they all decompose
in one-dimensional irreps. Further, the action of ~$\tB$~ is zero, thus, 
we speak only about the action of ~$\tD$. 
\item $\mu_A ~=~ 2n$, $n=1,2,...$ \\ 
For fixed ~$n$~ the representation decomposes into ~$2n+1$~
one-dimensional irreps. On one of these, which is spanned by the element:  
\eqn{triv} w_0 ~=~ \sum_{k=0}^n \binom{n}{k} \ta^{2n-2k} \td^{2k}
\ , \end{equation} 
the generator ~$\tD$~ acts by zero. The rest of the irreps 
are enumerated by the parameters:
~$\pm\ ,\ \tau$, where ~$\t = 2,4,...,2n=\mu_A\,$, and are
spanned by the vectors:
\eqnn{twod}&& u^\pm_\t ~=~ u_0 \pm \t u_1 \ ,\\ 
&& u_0 ~=~ \sum_{k=0}^n \a_k\, \ta^{2n-2k} \td^{2k}\ , \quad 
\a_0=1\ , \nn \\
&& u_1 ~=~ \sum_{k=0}^{n-1} \b_k\, \ta^{2n-2k-1} \td^{2k+1} \ , \quad 
\b_0=1\ , \nn
\eea 
on which ~$\tD$~ acts by:
\eqn{acttwo} 
\pi_R(\tD)\, u^\pm_\t ~=~ \pm\t \,u^\pm_\t 
 \end{equation}
which follows from: 
\eqn{acttwoo} 
\pi_R(\tD) \binom{u_0}{u_1} = \binom{\t^2\, u_1}{u_0} 
\end{equation}
Note that the value of the Casimir ~$\tD^2$~ is equal to $\t^2$.
The coefficients ~$\a_k\,,\b_k$~ depend on ~$\t$~ and 
are fixed from the two
recursive equations which follow from (\ref{acttwoo}):
\eqna{recu}
&&\t^2 \, \b_k ~=~ 2(n-k)\a_k - 2(k+1)\a_{k+1}\
, \quad k=0,...,n-1; \\ 
&&\a_k ~=~ (2k+1)\b_k - (2n-2k+1)\b_{k-1} 
\ , \quad k=0,...,n, 
\eena 
where we set ~$\b_{-1} \equiv 0$,\ ~$\b_n \equiv 0$. 
\item $\mu_A ~=~ 2n+1$, $n=1,2,...$ \\ 
For fixed ~$n$~ the representation is $(2n+2)$-dimensional and
decomposes into ~$2n+2$~ irreps which are
enumerated by two parameters: ~$\pm\ ,\ \t$, where 
~$\t = 1,3,5,...,2n+1=\mu_A\,$, and are spanned by the vectors: 
\eqna{twodd}&& w^\pm_\t ~=~ w_0 \pm \t w_1 \ ,\\ 
&& w_0 ~=~ \sum_{k=0}^n \a'_k\, \ta^{2n-2k+1} \td^{2k}\ , \quad 
\a'_0=1\ , \nn \\
&& w_1 ~=~ \sum_{k=0}^{n} \b'_k\, \ta^{2n-2k} \td^{2k+1} \ , \quad 
\b'_0=1\ , \nn
\eena 
on which ~$\tD$~ acts by (\ref{acttwo}). Note that the value of
the Casimir ~$\tD^2$~ is equal to $\t^2$. The coefficients
~$\a'_k\,,\b'_k$~ are fixed from the two recursive equations
which follow from (\ref{acttwo}): 
\eqnn{recuu}
&&\t^2\b'_k ~=~ (2n-2k+1)\a'_k - 2(k+1)\a'_{k+1} 
\ , \quad k=0,...,n; \\
&&\a'_k ~=~ (2k+1)\b'_k - 2(n-k+1)\b'_{k-1} 
\ , \quad k=0,...,n, 
\eea 
where we set ~$\a'_{n+1} \equiv 0$,\ ~$\b'_{-1} \equiv 0$. 
\end{itemize}

To summarize the list of irreps of ~$s14$~ on ~$S14$~ is:
\begin{itemize}
\item one-dimensional trivial 
\item two two-dimensional with all Casimirs ~$\tA, \tB^2,\tD^2$~
having the value 1.
\item one-dimensional enumerated by $n=1,2,...$ which for
fixed $n$ have Casimir values ~$2n,0,0$~ for ~$\tA,
\tB^2,\tD^2$, respectively.
\item pairs of one-dimensional irreps enumerated by $n=1,2,...$,
$\t=2,4,...,2n$, which have Casimir values ~$2n,0,\t^2$~ for
~$\tA, \tB^2,\tD^2$, respectively. 
\item pairs of one-dimensional irreps enumerated by $n=1,2,...$, 
$\t=1,3,...,(2n+1)$, which have  Casimir values ~$2n+1,0,\t^2$~
for ~$\tA, \tB^2,\tD^2$, respectively. 
\end{itemize}

Finally, we note in the irreps of ~$s14$~ on ~$S14$~ all
Casimirs can take only nonnegative integer values.

\section{Algebra $S14o$}
\label{sect:S14o}
\setcounter{equation}{0}

\subsection{Bialgebra relations}

In this Section we consider the matrix bialgebra ~$S14o$. 
We obtain it by applying the RTT relations (\ref{rtt})
for the case $\ R\ =\ R_{S1,4}\ $, cf. (\ref{eq:S14}), 
when ~$q^2= 1$. We shall consider the case ~$q=1$
(the case $q=-1$ is equivalent, cf. below). 
For $q= 1$ the relations following from (\ref{rtt}) and
(\ref{eq:S14}) are: 
 \eqn{eq:S14relq=1}
 a^2 = d^2  \ , ~~ b^2 = c^2 = 0 \ , ~~~~
ab = ba = ac = ca = bd = db = cd = dc = 0 
\end{equation}
or in terms of the generators $\ta,\tb,\tc,\td$~: 
\eqn{eq:S14trelq=1}
  \tb\ta = \ta\tb\ , ~~
 \tc\ta = - \ta\tc \ , ~~
  \td\ta = - \ta\td \ , ~~
 \tc\tb = - \tb\tc \ , ~~
  \td\tb = - \tb\td\ , ~~
 \td\tc = \tc\td 
\end{equation}
(The case $q=-1$ is obtained from the above through the exchange
~$\tb\leftrightarrow \tc$.)

{}From the above relations it is clear that we can choose any
ordering of the PBW basis. For definiteness we choose for the PBW
basis of $S14o$~:
\eqn{bas14o} \ta^k \tb^\ell \tc^m \td^n 
\end{equation}

\subsection{Dual algebra}

Let us denote by \ $s14o$\ the unknown yet dual algebra of $S14o$, 
and by ~$\tA,\tB,\tC,\tD$~ the four generators of \ $s14o$. 
We define the pairing \ $\lg Z, f\rg$, \
$Z=\tA,\tB,\tC,\tD$, \ $f$\ is from (\ref{bas14o}), as
(\ref{duda}). Explicitly, we obtain:
\eqna{eq:dualityS14o}
 &&\lg \ \tA \ , \ f \ \rg \ =\ 
 \ve\left(\frac{\partial f}{\partial \ta}\right) \ =\ 
 \begin{cases}k \qquad &{\rm for}\ \ f= \ta^k \\ 
 0 \qquad &{\rm otherwise }
 \end{cases} \\ 
 &&\lg \ \tB \ , \ f \ \rg \ =\ 
 \ve\left(\frac{\partial f}{\partial \tb}\right) \ =\ 
 \begin{cases}1 \qquad &{\rm for}\ \ f=\ta^k\tb \\ 
 0 \qquad &{\rm otherwise }
 \end{cases} \\ 
 &&\lg \ \tC \ , \ f \ \rg \ =\ 
 \ve\left(\frac{\partial f}{\partial \tc}\right) \ =\ 
 \begin{cases}1 \qquad &{\rm for}\ \ f= \ta^k \tc \\ 
 0 \qquad &{\rm otherwise }
 \end{cases} \\ 
 &&\lg \ \tD \ , \ f \ \rg \ =\ 
 \ve\left(\frac{\partial f}{\partial \td}\right) \ =\ 
 \begin{cases}
 1 \qquad &{\rm for}\ \ f=\ta^k\td \\ 
 0 \qquad &{\rm otherwise }
 \end{cases} 
\eena

Using the above we obtain:

\noindent 
{\bf Proposition 3:}\ The generators $\ \tA,\tB,\tC,\tD\ $
introduced above obey the following relations:
\eqna{eq:dualS14o}
 && [\tA,Z] = 0\ , \qquad \qquad Z=\tB,\tC,\tD\\
 && [\tB,\tC] = -2\tD \ , \quad [\tB,\tD] = -2\tC \ , \quad 
  [\tC,\tD] = -2\tB
\eena
\eqna{eq:dualcoprS14o}
 \delta_\cU(\tA) &=& \tA \otimes 1_\cU + 1_\cU\otimes \tA \\
 \delta_\cU(\tB) &=& \tB \otimes 1_\cU + 1_\cU \otimes \tB \\
 \delta_\cU(\tC) &=& \tC \otimes K + 1_\cU \otimes \tC 
\ , \qquad K=(-1)^{\tA}\\
 \delta_\cU(\tD) &=& \tD \otimes K + 1_\cU \otimes \tD 
\eena 
\eqn{coug} \ve_\cU (Z) = 0 \ , \qquad Z \ =\ \tA, \tB, \tC, \tD 
\end{equation} 
\eqn{anti}
 \g_\cU(\tA) = -\tA \ , \ \  
 \g_\cU(\tB) = -\tB \ , \ \ 
 \g_\cU(\tC) = - \tC K  \ , \ \ 
 \g_\cU(\tD) = - \tD K \ .
\end{equation} 
\PR The proof of (\ref{eq:dualS14o}) goes as the standard duality
between the classical $U(gl(2))$ and $GL(2)$, cf. e.g., \cite{Dod}. 
The proof of (\ref{eq:dualcoprS14o},\ref{coug},\ref{anti}) is also
standard, except the factor $K$ which appears while calculating: 
$$\lg \tC, \tc \ta^k \rg = (-1)^k \lg \tC, \ta^k\tc \rg= (-1)^k $$
which on the other hand is equal to (supposing an unknown yet $K$):
$$\lg \d_\cU (\tC), \tc\otimes \ta^k\rg = 
\lg \tC \otimes K + 1_\cU \otimes \tC \ ,\ \tc\otimes \ta^k\rg = 
\lg K, \ta^k \rg $$
Comparing the two RHSs we conclude that ~$K=(-1)^\tA$. The same
follows from calculating ~$\lg \tD, \td \ta^k \rg$.\spa 

\medskip 

\noindent 
{\bf Corollary 3:}\ The auxiliary generator $K=(-1)^\tA$ is central 
and ~$K^{-1} = K$. Its co-algebra relations are: 
\eqn{coalk} 
 \delta_\cU(K) = K \otimes K \ , \quad 
\ve_\cU (K) = 1\ , \quad \g_\cU(K) = K   \quad \spe 
\end{equation}

\medskip 

\noindent 
{\bf Corollary 4:}\ The algebra generated by the 
generator $\ \tA\ $ is a Hopf subalgebra of ~$s14o$. 
The algebra ~$s14o'$~ generated by $\ \tB,\tC,\tD\ $ 
is a subalgebra of ~$s14o$, but is not a Hopf subalgebra 
because of the operator $K$ in the co-algebra structure. 
The algebras ~$s14o,s14o'$~ are isomorphic to ~$U(gl(2)),U(sl(2))$,
respectively. The latter is seen from the following:
\eqna{sl2} 
&& X^\pm ~\equiv~ \half (\tD \mp \tC) \\ 
&& [\tB, X^\pm ] ~=~ \pm 2 X^\pm , \quad [X^+,X^-] ~=~ \tB \ . 
\eena 
Indeed the last line presents the standard $sl(2)$ commutation
relations. However, the generators ~$X^\pm$~ inherit the $K$
dependence in the coalgebra operations:
\eqna{antix}
&& \delta_\cU(X^\pm) ~=~ X^\pm \otimes K + 1_\cU \otimes X^\pm \\
&&\ve_\cU (X^\pm) ~=~ 0\\
&& \g_\cU(X^\pm) ~=~ - X^\pm K
\eena 
The algebra ~$s14o$~ is graded:
\eqn{grad}
\deg X^\pm\ =\ \pm 1\ , \qquad  
\deg \tA\ =\ \deg \tB\ =\ 0 \ , \qquad 
(\Lra\ \  \deg K\ =\ 0) \quad \spe 
\end{equation} 

\medskip 

Based on the above Corollary we are able to make the following important 
observation. The algebra ~$s14o$~ may be identified with a 
special case of the Hopf algebra ~$\cU_{p,q}$~ which was found 
in \cite{Dod} as the dual of $GL_{p,q}(2)$. To make direct contact  
with \cite{Dod} we need to replace there ~$(p^{1/2},q^{1/2}) \lra 
(p,q)$, then to set ~$q=p^{-1}$, and at the end to set ~$p =-1$. 
(The necessity to set values in such order is clear from, e.g., the 
formula for the co-product in (5.21) of \cite{Dod}.)  
The generators from  \cite{Dod} ~$K,p^K,H,X^\pm$~ correspond 
to ~$\tA,K,\tB,X^\pm$~ in our notation. 

More than this. It turns out that the corresponding 
algebras in duality, namely, ~$S14o$~ and ~$GL_{p,q}(2)$~ 
may be identified setting $q,p$ as above. To make this 
evident we make the following change of generators: 
\eqn{hatt} 
 \hat{a} \ =\  \ta+\tb, \quad
 \hat{b} \ =\  \td-\tc, \quad
 \hat{c} \ =\  \tc+\td, \quad
 \hat{d} \ =\  \ta-\tb\ .  
\end{equation} 
For these generators the commutation relations are: 
\eqn{eq:S14treq=1}
\hb\ha = -\ha\hb\ , ~\hc\ha = -\ha\hc\ ,
~\hd\ha =  \ha\hd \ , ~\hc\hb =  \hb\hc \ ,
~\hb\hd = -\hd\hb\ ,  ~\hc\hd = -\hd\hc 
\end{equation} 
i.e., exactly those of  ~$GL_{p,q}(2)$~ (cf. \cite{DMMZ}) 
for ~$p=q=-1$. Furthermore 
the co-product and and co-unit are as for $GL_{p,q}(2)$ or $GL(2)$, 
i.e., as in (\ref{coal}). For the antipode we have to suppose that 
the determinant ~$ad-p^{-1}bc$~ from \cite{DMMZ}, which here 
becomes (cf. $p=-1$): 
\eqn{omm} \omega = \ha\hd + \hb\hc \ ,\end{equation}
is invertible, or, that ~$\om\neq 0$~ and we extend the algebra 
by an element ~$\om^{-1}$~ so that:
\eqn{omma}  \om \om^{-1} ~=~ \om^{-1} \om ~=~ 1_\cA \ , \qquad 
\d_\cU(\om^{\pm 1}) ~=~ \om^{\pm 1} \otimes \om^{\pm 1} ~, 
~~~\ve_\cU(\om^{\pm 1}) ~=~ 1
 ~, ~~~\g_\cU(\om^{\pm 1}) ~=~ \om^{\mp 1}
\end{equation}
Then the antipode is given by:
\eqn{antip} 
\g_\cU \left( \begin{array}{cc}
\ha & \hb \cr \hc &\hd 
 \end{array}\right) ~=~ \om^{-1}\ 
\left(\begin{array}{cc}
\hd & \hb \cr \hc &\ha 
 \end{array}\right) 
\end{equation}
or in a more compact notation:
 \eqn{antipc} 
\g_\cU  \left(M\right) ~=~ M^{-1} 
\end{equation}
Indeed, we have: 
\eqn{antipp} 
\left( \begin{array}{cc} \ha & \hb \cr \hc &\hd  \end{array}\right) 
\left(\begin{array}{cc}\hd & \hb \cr \hc &\ha  \end{array}\right) 
~=~ 
\left(\begin{array}{cc}\hd & \hb \cr \hc &\ha  \end{array}\right) 
\left( \begin{array}{cc} \ha & \hb \cr \hc &\hd  \end{array}\right) 
~=~ \om\ 
\left(\begin{array}{cc}1 & 0 \cr 0 & 1  \end{array}\right) 
\end{equation} 

This relation between ~$s14o,S14o$~ and ~$\cU_{p,q}\,$, $GL_{p,q}(2)$~ 
was not anticipated since the corresponding ~$R$-matrices ~$R_{S1,4}$~ and 
~$R_{S2,1}$~ are listed in \cite{Hietarinta} as different and 
furthermore non-equivalent. It turns out that this is indeed the case, 
except in the case we have stumbled upon. To show this we first recall:
\eqn{s21} 
  R_{S2,1} = 
\left(
 \begin{array}{cccc}
 1 & 0 & 0 & 0 \cr
 0 & p & 1-pq & 0 \cr
 0 & 0 & q & 0 \cr
 0 & 0 & 0 & 1 \cr 
 \end{array}
 \right)
  \end{equation} 
which for ~$q=p^{-1}=-1$~ becomes:
\eqn{s21a} 
R_0 ~\equiv~ \left(  R_{S2,1}\right)_{q=p^{-1}=-1} 
~=~ \left(
 \begin{array}{cccc}
 1 & 0 & 0 & 0 \cr
 0 &-1 & 0 & 0 \cr
 0 & 0 &-1 & 0 \cr
 0 & 0 & 0 & 1 \cr 
 \end{array}
 \right)
  \end{equation} 
Further, we need:
\eqn{s14pm} 
R_\pm ~\equiv~ \left(R_{S1,4}\right)_{q=\pm 1} ~=~ 
\left(
 \begin{array}{cccc}
 0 & 0 & 0 & \pm 1 \cr
 0 & 0 & 1 & 0 \cr
 0 & 1 & 0 & 0 \cr
 \pm 1 & 0 & 0 & 0 \cr 
 \end{array}
 \right)
\end{equation} 

Now we can show that ~$R_\pm$~ can be transformed by  
"gauge transformations" to ~$R_0\,$, namely, we have:
\eqna{gaug} 
&& R_0 ~=~ (U_\pm \otimes U_\pm)\ R_\pm\ (U_\pm\otimes U_\pm)^{-1} \\[2mm]
&& U_+ ~=~ \frac{1}{\sqrt{2}} 
\begin{pmatrix} 1 & 1 \cr 1 & -1  \end{pmatrix} 
\ , \qquad 
 U_- ~=~ \frac{1}{\sqrt{2}} 
\begin{pmatrix} 1 & i \cr i & 1  \end{pmatrix} 
\eena 

In accord with this we have:
\eqnn{tut}
&&{\hat T} ~\equiv~ \begin{pmatrix} 
\ha & \hb \cr \hc & \hd  \end{pmatrix}\ , \qquad 
 T ~\equiv~ \begin{pmatrix} a & b \cr c & d  \end{pmatrix}
\ , \qquad 
{\hat T} ~=~ U_+\, T\, (U_+)^{-1} 
~~~~\Rightarrow \qquad \\
&& \ha ~=~ \half (a + b + c + d) \ , \qquad 
\hb ~=~ \half (a - b + c - d)\ , \label{tutu}\\
&& \hc ~=~ \half (a + b - c - d) \ , \qquad 
 \hd ~=~ \half (a - b - c + d)\nn  
\eea 
which is equivalent to substituting (\ref{eq:tatd}) 
in (\ref{hatt}).

The use of ~$U_-$~ would lead to different relations 
between hatted and unhatted generators, which, however, 
would not affect the algebra relations. Indeed: 
\eqnn{tutt}
&&{\hat T}' ~\equiv~ \begin{pmatrix} 
\ha' & \hb' \cr \hc' & \hd'  \end{pmatrix}\ , \qquad 
 T' ~\equiv~ \begin{pmatrix} a' & b' \cr c' & d'  \end{pmatrix}
\ , \qquad 
{\hat T}' ~=~ U_-\, T'\, (U_-)^{-1} ~~~~\Rightarrow \qquad \\
&& \ha' ~=~ \half (a' -i b' + i c' + d') \ , \qquad
 \hb' ~=~ \half (-i a' + b' + c' +i d') \ ,\\
&& \hc' ~=~ \half (i a' + b' + c' - i d') \ ,\qquad 
 \hd' ~=~ \half (a' +i b' -i c' + d') \nn 
\eea 
But this becomes equivalent to (\ref{tutu}) with the changes: 
\eqn{chhn}
 (\ha',i\hb',-i\hc',\hd') ~\mt~ (\ha,\hb,\hc,\hd) \ , \qquad
 (a',-i b',i c',d')  ~\mt~ (a,b,c,d) 
\end{equation}
while the (inverse) changes (\ref{chhn}) do not affect 
(\ref{eq:S14treq=1}),  (\ref{eq:S14relq=1}). 

\subsection{Representations of $s14o$ on $S14o$}
\label{sect:repr14o}

The regular representation of ~$s14o$, ($s14o'$)~ on itself and its weight
representations are the same as those of ~$U(gl(2))$, ($U(sl(2))$)~ due to
(\ref{sl2}). The situation is different for the 
representations of $s14o$ on $S14o$ since these involve the coalgebra 
structure. However, the deviation from the trivial coalgebra structure 
is only via the sign operator $K$, and as we shall see in some 
representations there remains no trace of this. 

In treating the 
representations of $s14o$ on $S14o$ we shall use the known 
construction for the induced representations of $\cU_{p,q}$ on $GL_{p,q}(2)$ 
from \cite{DoMi} and the relation between ~$s14o,S14o$~ and 
~$\cU_{p,q}\,$, $GL_{p,q}(2)$~ that we established in the previous subsection. 
For the comparison with \cite{DoMi} we should note the parametrization  
used there: ~$p ~=~ t s^{1/2}$, ~$q ~=~ t s^{-1/2}$. Thus, 
using ~$q=p^{-1}$, $p=-1$~ we need to substitute: ~$t\to 1$, ~$\sqrt{s} 
~\mt~ -1$. Further, one should substitute the operator ~$h$~ from 
\cite{DoMi} with ~$t^{\tB/2}$, and expand in order to get the action 
of ~$\tB$. Finally, one should substitute the operator ~$r$~ from 
\cite{DoMi} with ~$\sqrt{s}^\tA = (-1)^\tA = K$. In fact it is 
easier to derive the necessary formulae directly, which we shall 
proceed to do in a compact way.
 
Her we shall employ both the left action and the right action. 
We start by calculating the left  action using:
\eqna{lraa} 
&& \pi_L(Z) f ~\equiv~ \lg \g_\cU (Z),f_{(1)}\rg\,f_{(2)} 
\ , \qquad Z\in s14\,, \ Z\neq 1_\cU \,, \quad f\in S14 \ , \\
&&\pi_L(1_\cU) f ~\equiv~ f \ , \qquad f\in S14 \ .
\eena 
using for the PBW basis: 
\eqn{bas14a} \ha^j\, \hd^k\, \hb^\ell\, \hc^n\ . 
\end{equation}
For the left action on the elements of ~$S14o$~ we use a
Corollary of (\ref{lraa}):
\eqn{lrap} 
 \pi_L(Z) fg ~=~ \pi_L(\d'_\cU(Z)) (f\otimes g) 
\end{equation} 
where is used the ~{\it opposite comultiplication} 
~$\d'_\cU \equiv \s\circ\d_\cU\,$, where $\s$ is the permutation 
in ~$\cU \otimes\ \cU\,$. We find: 
\eqna{la14o} 
\pi_L(A)
\left(\begin{array}{ll} 
\ha^k & \hb^k\\ \hc^k & \hd^k 
\end{array} \right) 
&=& -k\,\left(\begin{array}{ll}
\ha^k & \hb^k\\ \hc^k & \hd^k
\end{array} \right) \\ 
\pi_L(K)\left(\begin{array}{ll} 
\ha^k & \hb^k\\ \hc^k & \hd^k
\end{array} \right) &=& (-1)^k \left(\begin{array}{ll} 
\ha^k & \hb^k\\ \hc^k & \hd^k
\end{array} \right) \\
\pi_L(B)\left(\begin{array}{ll} 
\ha^k & \hb^k\\ \hc^k & \hd^k 
\end{array} \right) &=& k\,\left(\begin{array}{cc}
-\ha^k &-\hb^k\\ \hc^k & \hd^k
\end{array} \right) \\ 
 \pi_L(X^+) \left(\begin{array}{ll} 
\ha^k & \hb^k\\ \hc^k & \hd^k 
\end{array} \right) &=& 
k\, \left(\begin{array}{cc} 
(-1)^{k-1} \ha^{k-1} \hc & \hd \hb^{k-1} \\ 0 & 0 
\end{array} \right) \\
\pi_L(X^-)\left(\begin{array}{ll}
\ha^k & \hb^k\\ \hc^k & \hd^k 
\end{array} \right) &=& 
k\,\left(\begin{array}{cc}
0 & 0 \\ \ha\hc^{k-1} &(-1)^{k-1}  \,\hd^{k-1} \hb   
\end{array} \right) \\
\pi_L(K)\, 1_\cA ~=~ 1 \ , && 
\pi_L(Z)\, 1_\cA ~=~ 0, \quad Z=\tA,\tB,\tC,\tD 
\eena  
(We give also the action of $K$ though it follows from that of $\tA$.) 

Next we calculate the right action as in (\ref{rraa}) 
to find:
\eqna{ra14o} 
 \pi_R(A)\left(\begin{array}{ll} 
\ha^k & \hb^k\\ \hc^k & \hd^k
\end{array} \right) &=& k\,\left(\begin{array}{ll} 
\ha^k & \hb^k\\ \hc^k & \hd^k
\end{array} \right) \\ 
 \pi_R(K)\left(\begin{array}{ll} 
\ha^k & \hb^k\\ \hc^k & \hd^k
\end{array} \right) &=& (-1)^k\left(\begin{array}{ll} 
\ha^k & \hb^k\\ \hc^k & \hd^k
\end{array} \right) \\
 \pi_R(B)\left(\begin{array}{ll} 
\ha^k & \hb^k\\ \hc^k & \hd^k 
\end{array} \right) &=& k\,\left(\begin{array}{ll} 
\ha^k & - \hb^k  
\\ \hc^k & -\hd^k  
\end{array} \right) \\
 \pi_R(X^+)\left(\begin{array}{ll} 
\ha^k & \hb^k\\ \hc^k & \hd^k 
\end{array} \right) &=&
k\, \left(\begin{array}{cc} 
0 & (-1)^{k-1} \ha \hb^{k-1}  \\ 
0 & \hd^{k-1} \hc   
\end{array} \right) \\ 
 \pi_R(X^-)\left(\begin{array}{ll} 
\ha^k & \hb^k\\ \hc^k & \hd^k 
\end{array} \right) &=&
k\,\left(\begin{array}{cc} 
\ha^{k-1}\hb & 0\\ 
(-1)^{k-1} \hd \hc^{k-1} & 0 
\end{array} \right) \\
\pi_R(K)\, 1_\cA ~=~ 1 \ , && 
\pi_R(Z)\, 1_\cA ~=~ 0, \quad Z=\tA,\tB,\tC,\tD 
\eena

By (\ref{lraa}) and (\ref{ra14o}) we have defined left and right action 
of ~$s14o$~ on $S14o$.  As in the classical case the left and right actions 
commute, and we shall use the right action to reduce the left
regular representation. Following \cite{Doa} we would like the right
action  to mimic some properties of a lowest weight module,
i.e., annihilation by the lowering (negative grade) generator 
$X^-$ and scalar action by the (exponent of the) 
Cartan (zero grade) generators ~$\tA$ (or $K$) and $\tB$. Such action 
is the reason we call these representations induced.  
We start with functions which are formal power series in the PBW basis:
\eqn{ser} f ~=~ \sum_{{j,k,\ell,m\,\in\,\ZZ_+}} 
\mu_{j,k,\ell ,m} ~\ha^j\, \hd^k\, \hb^\ell\, \hc^m 
\end{equation} 
The right-action conditions we mentioned are:
\eqna{rf} &&\pi_R(X^-)\ f ~=~ 0 \\
&&\pi_R(\tA)\ f ~=~ \r\,f \ ,\qquad  \pi_R(\tB)\ f ~=~ -\nu\,f
\eena
{}From (\ref{rf}a) follows that our functions would not depend on ~$\ha$~ 
and ~$\hc$, except through the determinant $\om$, 
since ~$\pi_R(X^-)\ \omega ~=~ 0$.  We also have: 
\eqn{raw} 
\pi_R(\tA)\,\omega^n = 2n \omega^n\ , \quad 
\pi_R(K)\,\omega^n =  \omega^n\ , \quad 
\pi_R(B)\,\omega^n = 0, \quad
\pi_R(X^\pm)\,\omega^n = 0 \ .
\end{equation} 
Thus, we continue with:
\eqn{serr} f ~=~ \sum_{{k,\ell,n \in\ZZ_+ \atop n\in\ZZ}} 
\mu_{k,\ell}\ \hd^k\, \hb^\ell\, \om^n
\end{equation} 
on which conditions (\ref{rf}b) lead to: ~$k+\ell+2n = \r\in\ZZ\,$, 
~$k+\ell=\nu\in\ZZ_+\,$, ~$\r-\nu\in 2\ZZ\,$, 
and the summation becomes single:
\eqn{serrr} f ~=~ \sum_{{\ell \in\ZZ_+}} 
\mu_{\ell}\ u_\ell \ , \qquad u_\ell ~\equiv~ 
\hb^\ell\, \hd^{\nu-\ell}\, \om^{(\r-\nu)/2} 
\end{equation} 
(where we changed the ordering since this would give simpler 
formulae for the action).  
Now if neither $\hb$ or $\hd$ has an inverse the representations will be 
finite-dimensional, in contrast to the classical case. 
However, these finite-dimensional representations we shall obtain also 
if we suppose that either $\hb$ or $\hd$ has an inverse (see below), 
and in the same time we shall have infinite-dimensional representations. 
Thus, further we shall suppose that ~$\hd$~ has an inverse. 
This means that we can allow ~$k$~ in (\ref{serrr}) to take any integer 
values, and then the same is true for ~$\nu$.

Now we shall work out the representation (left) action for the basis 
~$u_\ell$~ for which we need first the left action on ~$\om$~:
\eqn{law} 
\pi_L(A)\,\omega^n =-2n \omega^n\ , \quad 
\pi_L(K)\,\omega^n = \omega^n\ , \quad 
\pi_L(B)\,\omega^n = 0, \quad
\pi_L(X^\pm)\,\omega^n = 0 \ .
\end{equation} 
We also remark that the action on ~$\hd^{k}, \om^{n}$~ for negative $k,n$ 
is given again by (\ref{la14o}), (\ref{law}). (This can be checked, e.g.,  
by calculating ~$\pi_L(Z)\, \hd^{-k}\, \hd^{k} = \pi_L(Z)\,1_\cA$~ 
in two different ways: (\ref{lrap}) for the LHS, and 
(\ref{la14o}) for the RHS.) 
Then the  rules are: 
\eqna{transf}
\pi_{\nu,\r} (\tA)\, u_\ell &=& - \r\,u_\ell \\
\pi_{\nu,\r} (\tB)\, u_\ell &=& (\nu-2\ell) u_\ell \\
\pi_{\nu,\r} (X^+)\, u_\ell &=& (-1)^{\ell-1} \ell\, u_{\ell-1}\\
\pi_{\nu,\r} (X^-)\, u_\ell &=& (-1)^{\ell+1} (\ell-\nu)\, u_{\ell+1} 
\eena 
where ~$\pi_{\nu,\r}$~ denotes ~$\pi_L$~ with  
the parameter dependence made explicit. 

Thus, we have obtained infinite-dimensional representations of\ $s14o$\ 
parametrized by two integers ~$\nu,\r$. We shall denote by $\cC_{\nu,\r}$ 
the corresponding representation space. 
Note that these representations are highest weight representations 
since we have: ~$\pi_L (X^+)\, u_0 ~=~ 0$. If the parameter ~$\nu$~ 
is nonnegative, ~$\nu\in\ZZ_+$, then the corresponding representation 
is reducible, due to the fact that  ~$\pi_L (X^-)\, u_\nu ~=~ 0$. 
Thus, the vectors ~$u_0\,,...,u_\nu$~ form an invariant subspace, 
of dimension ~$\nu+1$, which shall denote  by ~$\cE_{\nu,\r}\,$, 
$\nu\in\ZZ_+\,$. Thus, if $\nu\in\ZZ_+$ we have two irreducible 
representations with representation spaces isomorphic to ~$\cE_{\nu,\r}$~ 
and to ~$\cC_{\nu,\r}/\cE_{\nu,\r}\,$ (the latter 
is infinite-dimensional). If $\nu\notin\ZZ_+$ the representation 
$\cC_{\nu,\r}$ is irreducible. 

{}From the above we are prompted to use the variable 
~$\y ~\equiv~ \hb \hd^{-1}$. This is also related to the 
following Gauss decomposition of $S14o$~:
\eqn{gag} 
\left(\begin{array}{ll} \ha & \hb \cr \hc &\hd
\end{array} \right) ~=~
\left(\begin{array}{cc} 
1 & \hb \hd^{-1}  \cr 0 &1
\end{array} \right)
\left(\begin{array}{cc} 
\om \hd^{-1}  & 0 \cr 0 &  \hd  
\end{array} \right)
\left(\begin{array}{cc} 
1 & 0 \cr \hd^{-1}\hc &1
\end{array} \right)
\end{equation} 
i.e., from here the natural variables are ~$\y$, $\hd$, $\om$.  
Thus, we use also the functions:
\eqn{serz} \vf ~=~ \sum_{{\ell \in\ZZ_+}} 
\a_{\ell}\ v_\ell \ , \qquad 
v_\ell ~\equiv~ \y^\ell\, \hd^{\nu}\, \om^{(\r-\nu)/2} 
\end{equation} 
The action of the generators on the variable ~$\y$~ is:
\eqn{lay} 
\left(\begin{array}{cc} 
\tA &  K\\ \tB   & X^\pm 
\end{array} \right) \y^\ell ~=~ 
\left(\begin{array}{cc} 
0 & \y^\ell \\ -2\ell\, \y^\ell & \pm \ell\,\y^{\ell\mp 1} 
\end{array} \right) 
\end{equation}
and the action on the basis is as for ~$u_\ell$~ except for ~$X^\pm$~: 
\eqna{transv}
\pi_{\nu,\r} (X^+)\, v_\ell &=&  \ell\, v_{\ell-1}\\
\pi_{\nu,\r} (X^-)\, v_\ell &=& (\nu-\ell)\, v_{\ell+1} 
\eena 
Thus, in this basis there is no trace of the non-triviality 
of the co-product of ~$X^\pm$. More than this we can reduce 
the representations directly to the classical ~$U(gl(2))$~ if we 
introduce the restricted functions ~$\hv(\y)$~ by the operators: 
\eqnn{res}&& \hA ~:~ \cC_{\nu,\r} ~\lra~ \hC_{\nu,\r}
~, \qquad \hv(\y) ~=~ \( \hA\,\vf \) (\y) ~\equiv 
~\vf (\y, 1_\cA, 1_\cA) \\
&&\hI ~:~ \hC_{\nu,\r} ~\lra~ \cC_{\nu,\r} ~, \qquad
\vf (\y,\hd,\om) ~=~ \( \hI\, \hv \) (\y,\hd,\om) ~\equiv ~ 
\hv(\y)\, \hd^\nu\, \om^{(\r-\nu)/2} \nn
\eea
We  denote the representation space of ~$\hv(\y)$~ by ~$\hC_{\nu,\r}$~ and
the representation acting in ~$\hC_{\nu,\r}$~ by ~$\hat\pi_{\nu,\r}$. 
The properties of 
~$\hC_{\nu,\r}$~ follow from the intertwining requirements \cite{Doa}: 
\eqn{int} \hpi_{\nu,\r}\ \circ\ \hA ~=~ 
\hA \ \circ\ \pi_{\nu,\r} ~, \qquad \pi_{\nu,\r}\ \circ\ \hI ~=~ 
\hI \ \circ\ \hpi_{\nu,\r} 
\end{equation}
In particular, the representation action of 
~$\hat\pi_{\nu,\r}$~ on ~$\eta^\ell$~ is given by the same formulae 
as the action of ~$\pi_{\nu,\r}$~ on ~$v_\ell\,$. 

At this moment, we should note that since we have functions 
of one variable $\eta$ we can treat it as 
complex variable $z$. In these terms we  recover 
from the action of ~$\hat\pi_{\nu,\r}$~ 
the classical vector-field representation of ~$gl(2)$~ 
(with $\pd_z \equiv d/dz$)~:
\eqn{hnn} \tA\,\hv ~=~ -\r\,\hv ~, \quad 
\tB\,\hv ~=~ (\nu - 2z\pd_z)\,\hv  ~, 
\quad X^+\,\hv ~=~  \pd_z\,\hv ~, \quad
X^-\,\hv = (\nu z - z^2\pd_z)\,\hv 
\end{equation}

Of course, the importance of the non-trivial co-product for ~$X^\pm$~ 
will be felt in the construction of the tensor products of the 
representations.

\section{Conclusions and outlook} 
\setcounter{equation}{0}

In this paper we have found the exotic matrix bialgebras which correspond 
to the two non-triangular nonsingular $4\times 4$ $R$-matrices of 
\cite{Hietarinta}, namely, ~$R_{S0,3}$~ and ~$R_{S1,4}$~ which 
are not deformations of the trivial $R$-matrix.
We study three bialgebras denoted by: ~$S03$, 
~$S14$, ~$S14o$, the latter two cases corresponding to  
~$R_{S1,4}$~ for deformation parameter ~$q^2\neq 1$~ and ~$q^2= 1$, 
respectively. We have found the corresponding dual bialgebras 
~$s03$, ~$s14$, ~$s14o$, and studied their 
representation theory. 

For the bialgebras ~$s03$~ and ~$s14$~ we have studied the 
regular representation (the algebra acting on itself), the 
weight representations, and the representations in which 
the algebra acts on the dual matrix bialgebra. The representation 
theory is degenerate: the irreps are finite-dimensional of 
maximal dimension 4 and 2 for $s03$ and $s14$, respectively. 
For future use we shall say that the bialgebras ~$s03,S03$~ 
and ~$s14,S14$ are ~{\it exotic}~ (adding to the list of 
exotic bialgebras termed so in \cite{ACDM1}).

The algebras ~$s14o,S14o$~ turned out to be Hopf algebras 
and to be  special cases of the two-parameter deformations 
~$\cU_{p,q}\,$, $GL_{p,q}(2)$, 
namely, they would be obtained from the latter by setting 
~$q=p^{-1}$~ and then $p=-1$. This was not anticipated since 
the corresponding $R$-matrices were different and seemingly 
nonequivalent cases of the classification of \cite{Hietarinta}. 
Thus, this became a case study  important methodologically, 
and so we have made the exposition according to the way 
we proceeded. In fact, the algebra ~$s14o$~ is equivalent 
even to ~$U(gl(2))$, and the only nontriviality is in the 
Hopf algebra structure. Thus, the regular and 
weight representations are as those of $U(gl(2))$. 
The induced representations of ~$s14o$~ on ~$S14o$~ could also 
be extracted from the equivalence with the two-parameter 
~$p,q$~ deformations but their consideration is also 
important methodologically. 

To conclude, we should stress that with this paper 
we finalize the explicit classification of 
the matrix bialgebras generated by four elements.
There are altogether ~{\it nine}~ such bialgebras, 
four of which are quantum groups and are deformations 
of the classical algebras of functions on $GL(2)$ and 
$GL(1|1)$ (two in each case), and the other five bialgebras, 
which we call exotic, are not such deformations. 

Further, we would like to study the spectral decomposition and
Baxterisation of these exotic algebras and associated noncommutative
geometries, cf. \cite{ACDM3}.

\end{document}